\documentclass[journal]{IEEEtran}
\newcommand{\beq}{\begin{equation}}
\newcommand{\eeq}{\end{equation}}
\newcommand{\bsq}{\begin{subequations}}
\newcommand{\esq}{\end{subequations}}
\newcommand{\bq}{\begin{eqnarray}}
\newcommand{\eq}{\end{eqnarray}}
\newcommand{\bqn}{\begin{eqnarray*}}
\newcommand{\eqn}{\end{eqnarray*}}

\usepackage{mathptmx}       
\usepackage{helvet}         
\usepackage{courier}        
\usepackage{type1cm}        
\usepackage{makeidx}         
\usepackage{graphicx}        

\usepackage{multicol}        
\usepackage[bottom]{footmisc}
\usepackage{ifthen}
\usepackage{subfigure}
\usepackage[usenames,dvipsnames]{color}

\makeindex             

\usepackage{graphics} 
\usepackage{graphicx} 
\usepackage{epsfig} 
\usepackage{epstopdf}
\usepackage{mathptmx} 
\usepackage{times} 
\usepackage[cmex10]{amsmath}

\usepackage{mathtools}  
\usepackage{amssymb}  
\usepackage{amsthm}

\usepackage{amsfonts}
\usepackage{cite}
\usepackage{verbatim}
\usepackage{comment}
\usepackage{algorithm}
\usepackage{algorithmic}
\usepackage{multirow}
\usepackage{cite}
\usepackage{stfloats}
\usepackage{supertabular}
\usepackage{longtable}

\DeclareGraphicsExtensions{.pdf,.png,.jpg,.eps,.ps}
\renewcommand{\arraystretch}{1.2}
\usepackage{arydshln}
\usepackage{multicol}
\usepackage{colortbl}
\theoremstyle{definition}

\newtheorem{proposition}{Proposition}

\theoremstyle{definition}
\newtheorem{definition}{Definition}

\ifCLASSINFOpdf
\else
\fi

\usepackage{comment}
\usepackage{ifthen}
\newboolean{showcomments}
\setboolean{showcomments}{true}
\newcommand{\ychen}[1]{\ifthenelse{\boolean{showcomments}}
        { \textcolor{red}{YC: #1}}}
\newcommand{\fliu}[1]{\ifthenelse{\boolean{showcomments}}
        { \textcolor{blue}{(FL:  #1)}}{}}
\newcommand{\slow}[1]{\ifthenelse{\boolean{showcomments}}
        { \textcolor{red}{(SL: #1)}}}

\usepackage{amsmath}
\allowdisplaybreaks[4]

\renewcommand{\baselinestretch}{0.99}

\begin{document}

%
\title{An Energy Sharing  Game with Generalized Demand Bidding: Model and Properties}
%
%
%

\author{Yue Chen,
        Shengwei Mei,~\IEEEmembership{Fellow,~IEEE,}
        Fengyu Zhou,
        Steven H. Low,~\IEEEmembership{Fellow,~IEEE,}
        Wei Wei, ~\IEEEmembership{Senior Member,~IEEE,}
        and Feng Liu,~\IEEEmembership{Senior Member,~IEEE}
}

%
%

\markboth{Journal of \LaTeX\ Class Files,~Vol.~XX, No.~X, Feb.~2019}%
{Shell \MakeLowercase{\textit{et al.}}: Bare Demo of IEEEtran.cls for IEEE Journals}
%



\maketitle

\begin{abstract}
	This paper proposes a novel energy sharing mechanism for prosumers who can produce and consume. Different from most existing works, the role of individual prosumer as a seller or buyer in our model is endogenously determined. Several desirable properties of the proposed mechanism are proved based on a generalized game-theoretic model. We show that the Nash equilibrium exists and is the unique solution of an equivalent convex optimization problem. The sharing price at the Nash equilibrium equals to the average marginal disutility of all prosumers. We also prove that every prosumer has the incentive to participate in the sharing market, and prosumers' total cost decreases with increasing absolute value of price sensitivity. Furthermore, the Nash equilibrium approaches the social optimal as the number of prosumers grows, and competition can improve social welfare.  
\end{abstract}


%
\IEEEpeerreviewmaketitle

\section*{Nomenclature}
\addcontentsline{toc}{section}{Nomenclature}
\subsection{Indices and Sets}
\begin{IEEEdescription}[\IEEEusemathlabelsep\IEEEsetlabelwidth{sssssss}]
	\item[$i/n$]         Index of prosumers/resources.
	\item[$k$] Index of resources of a prosumer.
	\item[$\mathcal{I}/\mathcal{N}$] Set of prosumers/resources.
	\item[$\mathcal{S}/\mathcal{D}$] Set of sellers/buyers.
	\item[$\mathcal{K}_i$] Set of resources of prosumer $i$.
	\item[$f_i(\cdot)$] Disutility function of prosumer $i$.
	\item[$s_i(\cdot)$] Sharing cost function of prosumer $i$.
	\item[$md_i(\cdot)$] Marginal disutility of prosumer $i$.
	\item[$\Pi_i(\cdot)$] Total cost function of prosumer $i$.
	\item[$X_i$] Action set of player $i$, and $X = \prod_i X_i$.
	\item[$\Gamma_i(\cdot)$] Cost function of the sharing game.
\end{IEEEdescription}

\subsection{Parameters}
\begin{IEEEdescription}[\IEEEusemathlabelsep\IEEEsetlabelwidth{sssssss}]
	\item[$I/N$]    Number of prosumers/resources.
	\item[$D_i^0$] Fixed amount of energy prosumer $i$ consumes.
	\item[$p_i^0$] The amount of energy $i$ generates originally.
	\item[$p_i^{r0}$] Forecast output of random generation device $i$.
	\item[$\Delta p_i^{r0}$] Deviation of random generation device $i$.
	\item[$E_i^0$] The amount of energy bought from the grid by $i$.
	\item[$D_i$] Required load reduction for prosumer $i$.
	\item[$c_i,d_i$] Disutility function coefficients for prosumer $i$.
	\item[$c_i^k,d_i^k$] Disutility coefficients of resource $k$ of prosumer $i$.
	\item[$\underline D,\overline D$] Lower/upper bound of $D_i$.
	\item[$\underline c,\overline c$] Lower/upper bound of $c_i$.
	\item[$\underline d,\overline d$] Lower/upper bound of $d_i$.
	\item[$a$] Price sensitivity of prosumers.
	\item[$K_i$] Number of resources of prosumer $i$.
	\item[$\alpha$] Proportion of total cost reduction paid to platform.  
\end{IEEEdescription}

\subsection{Decision Variables}
\begin{IEEEdescription}[\IEEEusemathlabelsep\IEEEsetlabelwidth{sssssss}]
	\item[$p_i$] Output adjustment of prosumer $i$.
	\item[${\lambda _c}$] Sharing market clearing price.
	\item[$b_i$] Willingness to buy/pay of prosumer $i$.
	\item[$\bar{b}$] Average purchase desire of all prosumers.
	\item[$q_i$] Amount of energy bought/sell from/to the market.
	\item[$p_i^k$] Output adjustment of resource $k$ of  prosumer $i$.
	\item[$\xi$,$\xi^{'}$] Dual variable of the energy balance equation of the equivalent central decision-making problem.
\end{IEEEdescription}

\subsection{Abbreviations}
\begin{IEEEdescription}[\IEEEusemathlabelsep\IEEEsetlabelwidth{sssssss}]
	\item[DER] Demand-side energy resource.
	\item[SDF] Supply-demand function.
	\item[GNG] Generalized Nash game.
	\item[NE] Nash equilibrium.
	\item[MRP] Multi-resource prosumer.
	\item[IDL] Individual.
	\item[SMK] Sharing market.
	\item[SCO] Social optimal.
	\item[PoA] Price of anarchy.
\end{IEEEdescription}

\section{Introduction}
%
%
%
%
\IEEEPARstart{T}{he} proliferation of distributed wind and solar power has endowed traditional consumers with the capability of generation, precipitating the advent of ``prosumers'' \cite{prosumer}. Local trading among prosumers may provide the system with additional flexibility \cite{prosumer2}. Usually, prosumers are managed in a centralized manner. As the number of prosumers continues to increase, management inefficiency and potential conflict of interests may move the prosumer management towards a distributed and scalable way. In such a circumstance, the mismatch of supply and demand for individual prosumer turns out to be a core issue and the surplus energy trading among prosumers becomes a must \cite{zafar2018prosumer}. Different from traditional consumers, prosumers can choose either to buy or sell energy, motivating new market mechanisms for supporting flexible energy exchanging among prosumers.

Nowadays, online platforms and applications have enabled resource sharing in more and more sectors \cite{Sharing-develop1}, such as ride-sharing (e.g., Uber, Lyft), room-sharing (e.g., AirBnB), and workplace-sharing (e.g., Upwork). These sharing platforms allow people to provide their idle goods to someone in need for a profit. Existing studies in economics have investigated the operation of product sharing \cite{Sharing-Benifit2}. These successes motivate a new paradigm of energy utilization in power systems, where energy prosumers share their energy in a similar way \cite{Sharing-Benifit1}. The main difficulty of sharing platform construction is the design of an appropriate sharing mechanism, which means how to provide the energy, how to clear the market, and how to allocate the revenue. Performance of typical sharing platforms is studied in \cite{Sharing-platform4}, and the influence of prices and subsidies in \cite{Sharing-price}. A review of sharing economy can be found in \cite{Sharing-overview}.

For energy sharing, the potential of game-theoretic approaches was summarized in \cite{Sharing-analytical0}, including the applications in electric vehicles (EV), demand-side energy resource (DER) and storage managements. The economic efficiencies of individual, sharing and aggregation schemes were quantitatively compared in \cite{Sharing-analytical1}, showing that energy sharing can achieve near-optimal efficiency without a central coordinator, which is a promising scheme for future energy market organization. Random sharing clearing price in a storage investment problem was characterized in \cite{Sharing-analytical2} using a simplified time-of-use (TOU) model. Above work initially explores the problem and opportunity of sharing in smart grid, and the models are relatively abstract and simple. More detailed analytical studies related to resource sharing can be roughly cast into the following three categories. Here, ``side'' refers to the role of market participant (as a seller or a buyer) as in \cite{gerding2013two,zhang2016energy}.

\textbf{Two-sided market with clearing price.} The sellers report the amount of products to share or their cost coefficients; the buyers report the amount of products they want or the money they are willing to pay. With these bids, the third-party sharing platform solves an optimization problem with the objective of social welfare or self-revenue maximization and clears the market. Ref. \cite{Sharing-twosided1} provided insights into the tradeoff between revenue maximization and social welfare maximization. The connection between clearing price and the Vickrey-Clarke-Groves (VCG) mechanism was analyzed in \cite{Sharing-twosided2}. Incentive design for electric vehicle-to-vehicle charge sharing was investigated in \cite{Sharing-twosided3}. System constraints can be taken into account in the two-sided market analysis. However, since the supply and demand statuses of participants are predetermined, it cannot fully capture the behaviors of prosumers who can choose to buy or sell changeably.

\textbf{Single-sided market with set price.} Different from the two-sided market, it assumes that the roles of all participants are symmetric. The benefits from sharing are distributed among prosumers via prices set by the sharing platform. An hour-ahead optimal sharing pricing model for photovoltaic (PV) prosumers was proposed based on a Stackelberg game \cite{Sharing-setprice1}. Uncertainty was further taken into account in \cite{Sharing-setprice2} with two kinds of sharing schemes, i.e., the direct sharing (within one time period) and the buffered sharing (across different time periods). Stackelberg model was also used in \cite{el2017managing} to study the bounded rationality of prosumers. The results under traditional expected utility theory and prospect theory were analyzed and compared. Another research focus was the peer-to-peer (P2P) energy sharing \cite{Sharing-setprice3, paudel2018peer, long2018peer}. In the above studies, the sharing prices are set by the platform via solving a Stackelberg game, in which the upper level is the platform's pricing problem and the lower level prosumers' decision making problems. The impact of one prosumer's strategy on the other prosumers' decisions is not fully captured.

\textbf{Single-sided market with re-allocation.} In this kind of sharing, the benefit distribution is achieved via re-allocation instead of price regulation. The main difficulty stems from the design of re-allocation schemes. The renowned VCG mechanism \cite{Sharing-reallocation1} could be regarded as an example. Although the VCG re-allocation approach is ease to implement, it is not self-budget balancing as extra bonuses outside the sharing market is required. A cost re-allocation method for a group of electricity storages was presented in \cite{Sharing-reallocation2}, resulting in a cooperative game. A coalitional game based algorithm was proposed in \cite{mei2019coalitional} for energy exchange among microgrids. A conceptual design for the DERs sharing was proposed in \cite{Sharing-reallocation3}, where an aggregator coordinates all DERs in real-time operation and evaluates coordination surplus, which is split between aggregators and prosumers. However, the redistribution after a sharing transaction is difficult in practice, as it requires some private information of individual participants, e.g. the storage capacity or the cost coefficients.

\begin{figure}[!t]
	\centerline{\includegraphics[width=0.65\columnwidth]{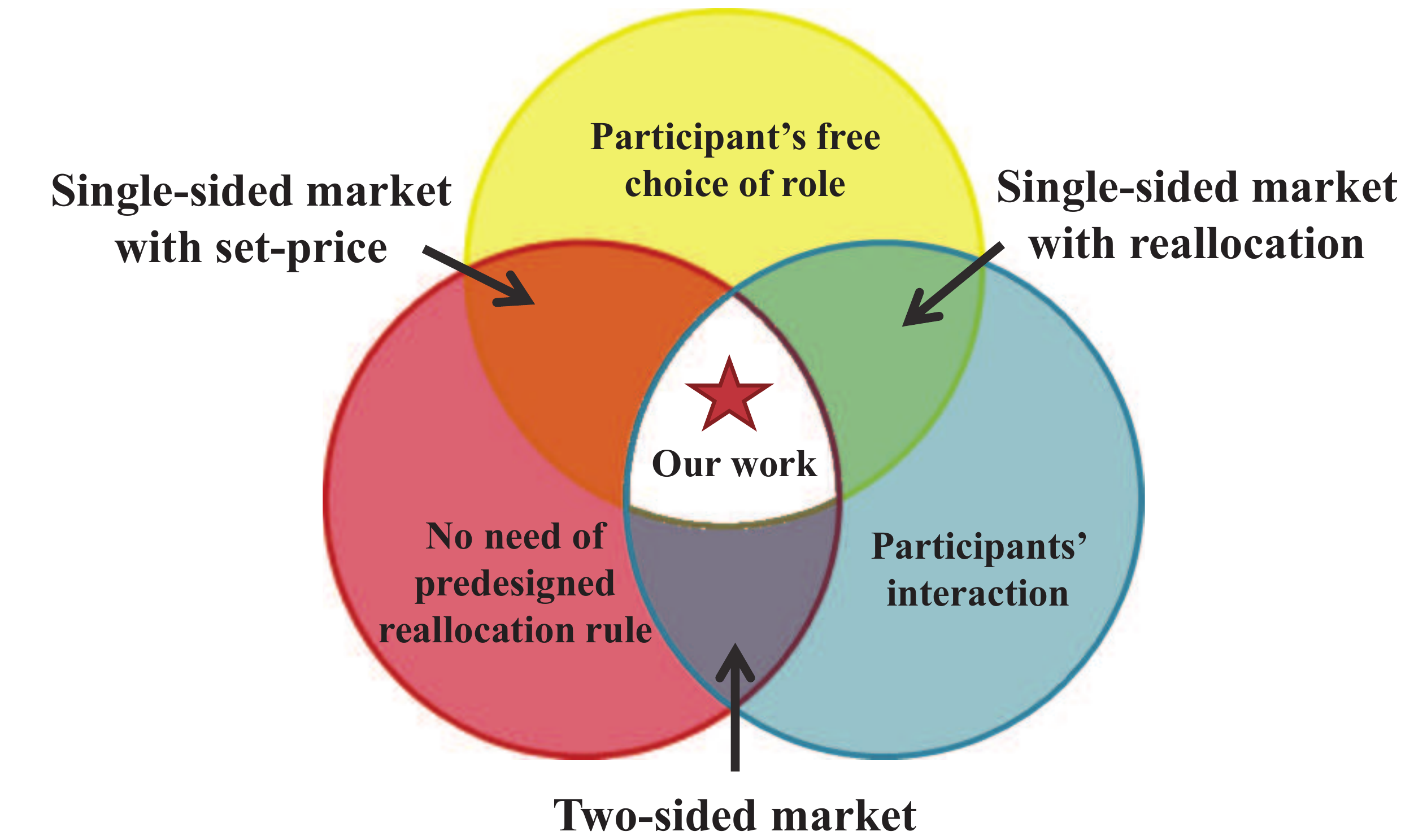}}
	\setlength{\abovecaptionskip}{0pt}
	\caption{Comparison of our work with current studies.}
	\setlength{\belowcaptionskip}{-5em}
	\label{fig:comparison}
	\vspace{-1.5em}
\end{figure}

Taking into account the unique features of the energy system, the requirements for an energy sharing market are as follows. The advantages of our work over existing studies are compared from three aspects as in Fig. \ref{fig:comparison}.

1) \emph{Each prosumer is free to sell or buy}. Compared with the ``two-sided market'', where the roles of participants are predetermined, the prosumers in our model are symmetric, which is in more accordance with the reality. Whether a prosumer will become a seller or a buyer is endogenously determined. To be specific, each prosumer only needs to bid a $b_i$, showing his willingness to buy. It has been proved that when a prosumer has less willingness to buy, he naturally becomes a seller; otherwise, he is a buyer.

2) \emph{The interaction among prosumers is fully activated.} With more and more prosumers, the online platform has been changing its role from a control center to a supporting platform. In the ``single-sided market with set-price'', the platform decides the price and prosumers act only as price-takers, which discourages the enthusiasm of prosumers to participate in sharing. In our model, each prosumer makes strategic bidding and the platform clears the market according to these bids, which creates sufficient incentives to encourage the participation of prosumers in sharing.

3) \emph{The sharing mechanism is easy to implement with good properties}. Compared with the ``single-sided market with reallocation'', our mechanism does not require private information and the market clearing rule is simple and transparent. Also, it possesses several provable properties, such as the existence and uniqueness of Nash equilibrium as well as the achievement of a Pareto improvement, to name a few.

This paper proposes a simple energy sharing mechanism based on generalized demand bidding. It is similar to the \emph{supply function bidding} \cite{Supply-function-bid1,Supply-function-bid2} in demand response programs. 
Under this mechanism, each seller submits his supply function to the auctioneer, then a market clearing price is set according to the submitted supply functions and the expected total load shedding. Supply function bidding can fully capture the impact of a seller's bid on his contracted quantity as well as the clearing price, and is effective in competitive markets \cite{Supply-function-bid3}. A mathematical problem with equilibrium constraints (MPEC) was used to depict the strategic bidding of a single-firm in \cite{hobbs2000strategic} and an iterative algorithm was adopted to obtain the equilibrium with multiple firms. The demand response program is modeled as two noncooperative games, one among suppliers based on supply function bidding mechanism, and the other among customers \cite{kamyab2015demand}. When it comes to the sharing market, the problem is more complicated in the following two aspects which are included in our framework: 1) instead of selling only, the prosumers weigh between self-consumption and sharing; 2) sellers and buyers, as well as other equilibrium quantities, are endogenously determined and their roles can change over time. 
In this regard, we generalize the supply function to a generic supply-demand function (SDF), based on which we build a sharing mechanism for energy prosumers. 
This work possesses three salient features:

1) \textbf{New mechanism for energy sharing}. An effective sharing mechanism based on a generic supply-demand function is proposed. In contrast to the ``two-sided market'' based analysis, each participant is free to act as a seller or a buyer, aiming to minimize his own disutility. Instead of setting price directly as in the ``single-sided market with set price'', the platform supports and balances mutual impacts among prosumers. The model that encapsulates the decision making of prosumers turns out to be a generalized Nash game (GNG), which is further reduced to a standard Nash game. The existence of a unique Nash equilibrium is proved.

2) \textbf{Provable properties of the sharing mechanism}. Properties of the sharing equilibrium price are presented. It is proved that a Pareto improvement can be achieved, implying that every prosumer's cost under sharing is never worse than that under individual decision-making. Thus, all prosumers have motivation to take part in sharing. Moreover, the total cost of all prosumers decreases with the increase of price sensitivity, and the sharing market will eventually lead to the same outcome as the social optimum with more and more prosumers. It is also revealed that the proposed sharing mechanism is budget self-balancing and no private information is needed for re-allocation. This makes it easier to implement compared with the re-allocation based schemes.

3) \textbf{Impacts of competition on social efficiency}. To simplify the analyses, the basic model of our mechanism is based on a perfectly monopolistic competition situation, in which every prosumer controls only one resource. We further investigate a more realistic case, in which a prosumer could possess multiple resources. A special case provides a proof of concept that spreading the resources among more prosumers subject to some appropriate rules may further reduce the social cost.

The rest of this paper is organized as follows. The mathematical formulation of energy prosumers and description of the energy sharing mechanism are presented in Section II. Some basic properties of the sharing game are given in Section III. The impact of competition on social welfare is studied in Section IV. Some possible extensions are provided in Section V and illustrative examples in Section VI. Finally, conclusions are summarized in Section VII.

\section{Game  Model of Energy Sharing}
\subsection{Energy Prosumers}
In this paper, we consider the decision-making problems of a set of prosumers $\mathcal{I}$,  indexed by $i \in \mathcal{I}=\{1,2...,I\}$. 
There are $N$ kinds of resources, indexed by $n \in \mathcal{N}=\{1,2,...,N\}$, which can be a distributed generator (DG), virtual power plant (VPP), and etc. First, we consider the case under perfectly monopolistic competition market, where each prosumer $i$ owns one kind of resource, and we have $I=N$. To distinguish from the case under imperfect monopolistic competition market, we use $N$ to represent the number of prosumers here. The fixed amount of energy prosumer $i$ consumes is $D_i^0$, and is satisfied by the energy it generates $p_i^0$, as well as the energy bought from the grid $E_i^0$. These prosumers take part in a demand response program, and the required load reduction for prosumer $i$ is a given value $D_i$ ($\underline D \le D_i \le \overline D$), which means the amount of energy it bought from the grid needs to be reduced by $D_i$. Each prosumer changes its resource output to meet the load adjustment. For example, to reduce its load by $D_i>0$, prosumer $i$ needs to increase output by $D_i$. Any deviation from the original operating point will cause disutility. The disutility function of prosumer $i$ is a quadratic function $f_i(p_i)=c_i p_i^2+d_i p_i$, where $p_i$ is the output adjustment of resource, and $0 < \underline c \le c_i \le \overline c$, $0 < \underline d \le d_i \le \overline d $ are the cost coefficients. The quadratic utility function is widely adopted in smart grid analysis \cite{samadi2012advanced}. Here we choose a quadratic function because it has some properties that well fit the characteristic of prosumer's disutility, i.e., both increasing and decreasing the output result in disutility; the prosumer's marginal disutility grows with increasing adjustment.

Generally, a demand response program of a utility company roughly consists of two main decisions: the first is how much of the total demand to adjust and how to allocate the total adjustment to prosumers in its service territory, and the second is how the prosumers fulfill the allocated load adjustment. The typical non-price based demand response addresses the first problem, while our model focuses on the second one. To be specific, first, the load adjustment requirement is sent to the utility from the bulk grid. Sufficient/insufficient load shed can be rewarded/penalized based on a contract. The utility operator solves an optimization problem to determine how much each prosumer should adjust considering these financial terms. Then, given the load adjustment requirement $D_i$, traditionally, each prosumer simply fulfill this command by reducing its own load by $D_i$ as requested. In this paper, we ask the question, if there is a platform that would allow (but not force) the prosumers to trade, can they fulfill their required reductions $D_i, \forall i$ in a way that is more profitable than without the sharing market, and if so, how?

When a prosumer $i\in \mathcal{I}$ takes part in a demand response program individually, there is no room for optimization since $p_i=D_i$ is clearly the solution. The corresponding cost is $f_i(D_i)$. However, this result may not be the most efficient if prosumers with different marginal disutilities are allowed to trade with others. In this context, the design of an effective profit allocation scheme, from which all prosumers take part in sharing can benefit, is desired. The traditional supply function bidding in demand response program cannot be applied because of the simultaneous non-deterministic clearing quantity and clearing price as a prosumer can changeably act as either a producer or a consumer. Hence a more general bidding mechanism, which can reflect prosumers' willingness to buy or sell energy while determining both the clearing quantity and price, are necessary.

\subsection{Generic Supply-Demand Function}
In this subsection, we propose a generic supply-demand function by generalizing the conventional supply function, so as to consider the situation where the participant can flexibly change his role between a seller and a buyer. 

In the sharing market, the demand (or supply) function \cite{hobbs2000strategic} of each prosumer can be expressed by
\bq
\label{eq:clearcondition1}
q_i=a_i\lambda_c+b_i 
\eq
where $\lambda_c$ is the market clearing price, $q_i$ is the amount of energy ($q_i>0$ means he is a buyer and gets energy from the sharing market, $q_i<0$ means he is a seller and sells energy to the sharing market). $a_i<0$ represents price sensitivity and $b_i$ shows his willingness to buy. For simplification, we assume all prosumers have the same price sensitivity, i.e. $a_i=a <0, i \in \mathcal{N}$. The average purchase desire is defined as $\bar{b}=(\sum_i b_i)/N$. The market clears when the net quantity $\sum \nolimits_i q_i=0$ and the obtained sharing price is
\bq
\label{eq:clearcondition}
\lambda_c=-{\sum \nolimits_i b_i}/{(Na)}=-{\bar{b}}/{a}
\eq
Here, $b_i \ge \bar{b}$ implies prosumer $i$ is more willing to buy than the average. We have $q_i=a\lambda_c+b_i \ge 0$, and the prosumer appears to be a buyer. Similarly, a prosumer who has less willingness to buy than the average ($b_i \le \bar{b}$) turns to be a seller ($q_i \le 0$). Therefore, the statuses of prosumers are determined spontaneously by their purchase desires, which enables a simple but effective sharing mechanism, as we explain.

\textbf{Remark 1:} The assumption that universal price sensitivity is used can be justified from the following perspectives.

1) \textbf{Physical and Implementation Interpretation}. Here, we conjecture that prosumers within a utility usually have similar price sensitivities, so it is acceptable to assume a universal price sensitivity for simplification. Even with such an assumption, the personal character of each prosumer is still retained since the prosumer can accordingly adjust his bid $b_i$ to show his willingness to buy. In fact, the $a$ can be regarded as the price sensitivity of the whole sharing market, reflecting how the prosumers' total bid $\sum_{i \in \mathcal{I}} b_i$  influences the sharing price. Moreover, when it comes to application, different $a_i$ is a private information and hard to accurately obtain. The universal $a$ can be regarded as a parameter in the set rule for market clearing and as proved later in Proposition \ref{Thm:prop-3}, even though this universal $a$ does not match each prosumer's price sensitivity accurately, he is still willing to take part in sharing.

2) \textbf{Model and Properties Extension}. Admittedly, it'd be a less restrictive model if we allow $a_i$ to be heterogeneous. Though we might not be able to prove some special properties of equilibrium, the sharing mechanism itself extends to heterogeneous $a_i$. To be specific, we can allow different but fixed $a_i,\forall i$, with each prosumer still only bidding his $b_i$. The clearing condition in \eqref{eq:clearcondition} becomes \eqref{eq:clearcondition2} and the sharing game as well as Proposition \ref{Thm:prop-2}-\ref{Thm:prop-1} in this paper extends accordingly.
\bq
\label{eq:clearcondition2}
\lambda_c = -\sum \nolimits_i b_i / \sum \nolimits_i a_i
\eq
For Proposition \ref{Thm:prop-4}, We run some simulations with heterogeneous $a_i$ in Section VI.B to show that it approximately holds.

\subsection{Energy Sharing Mechanism}
\begin{figure}[!t]
	\centerline{\includegraphics[width=0.65\columnwidth]{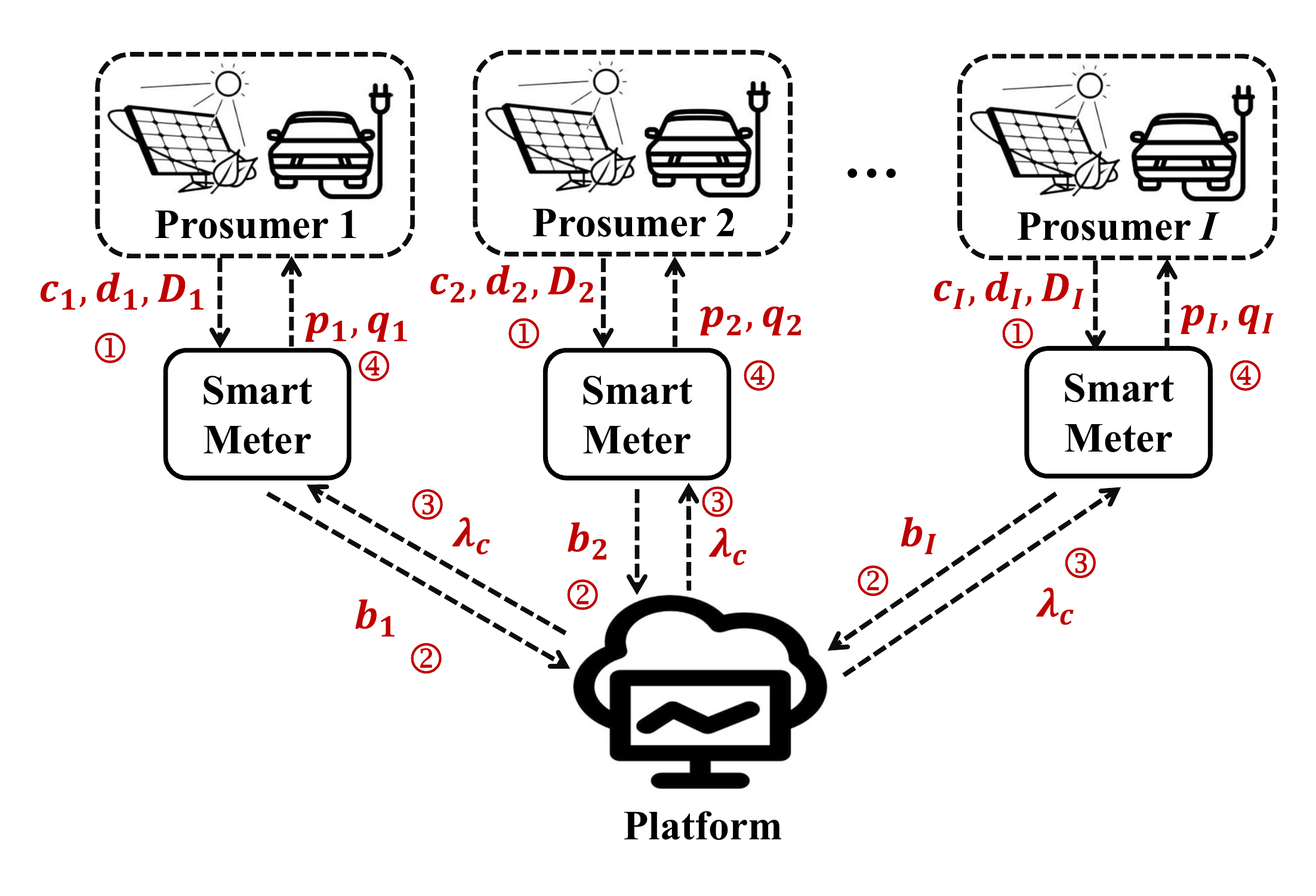}}
	\setlength{\abovecaptionskip}{0pt}
	\caption{Framework of the proposed sharing mechanism.}
	\setlength{\belowcaptionskip}{-5em}
	\label{fig:framework}
	\vspace{-1.5em}
\end{figure}
There are $N$ prosumers in the sharing game on a supporting platform. Each prosumer is connected to the platform via a smart meter with a bidirectional information channel as in Fig.\ref{fig:framework}. The information flow of the sharing game is as follows.

\textbf{Step 1}: Each prosumer inputs his cost coefficients $c_i,d_i$, and the required load reduction $D_i$ to the smart meter. These data are kept private and will not be accessed by the platform or other prosumers. Estimate the value of price sensitivity $a$ via historical data. Initialize $\lambda_c=0$.

\textbf{Step 2}: The smart meter $i$ decides a bid $b_i$ via the embedded optimization problem \eqref{eq:Nash game} with $\lambda_c$ (from which it can deduce $\sum_{j \ne i} b_j$) and sends it to the platform.

\textbf{Step 3}: After receiving all the bids $b_i,\forall i$, the platform clears the market by setting price to $\lambda_c(b) = -\sum_i b_i /(Na)$ and returns the sharing price $\lambda_c$ to each smart meter.

\textbf{Step 4}: If the sharing price $\lambda_c$ remains the same as in the last iteration, go to \textbf{Step 5}; otherwise, back to \textbf{Step 2} with the updated $\lambda_c$.

\textbf{Step 5}: The average purchase desire is $
\bar{b} = \sum \nolimits_i b_i / N  \nonumber
$. The smart meter calculates the self-production $p_i$ and sharing amount $q_i(b)$, sends it back to prosumer $i$ and the prosumer executes it. If $b_i \ge \bar{b}$, $q_i(b) \ge 0$, which means prosumer $i$ will buy $q_i(b)$ from the sharing market and his payment is $\lambda_c(b)q_i(b)$. Otherwise, if $b_i \le \bar{b}$, $q_i(b) \le  0$, which means prosumer $i$ will sell $-q_i(b)$ to the sharing market and he will get $-\lambda_c(b)q_i(b)$. It is worth noting that, since $q_i(b) \le 0$, $-q_i(b)$ is a positive quantity and selling $-q_i(b)$ indicates the prosumer is a seller.

Under this setting, as shown in Fig.\ref{fig:framework}, prosumer $i$ only needs to provide his private information ($c_i, d_i, D_i$) to his own smart meter, which will not be known by other prosumers or the platform, and so privacy is well preserved.
Besides, Instead of a one-time bidding, an iterative process is adopted where all smart meters bid $b_i,\forall i$ to the platform, then the platform clears the market and sends $\lambda_c$ back to each smart meter, and then the smart meters adjust their bids, and so on, until the equilibrium is reached.

The sharing market clears when
\bq
\label{eq:clearing}
\sum\nolimits_{i \in \mathcal{S}} (-q_i)=\sum \nolimits_{i \in \mathcal{D}} q_i
\eq
$\mathcal{S}$ is the set of sellers, $\mathcal{D}$ is the set of buyers. Each prosumer $i$ belongs to either $\mathcal{S}$ or $\mathcal{D}$, which means $\mathcal{I}=\mathcal{S}\cup\mathcal{D}$. Hence equation \eqref{eq:clearing} also implies
\bq
\label{eq:clearing-2}
\sum \nolimits_{i \in \mathcal{I}} (a\lambda_c+b_i)=0
\eq

The proposed sharing mechanism is essentially a kind of auction approach. Other examples of auction approach in economy are English auction, Dutch auction and first-price sealed-bid auction, to name a few. However, these auction approaches are designed for many participants competing for one thing, while in this paper, each prosumer gets some amount of energy.
\vspace{-1em}
\subsection{Energy Sharing as A Generalized Nash Game}
It is easy to verify that the setting price $\lambda_c(b)$ clears the market.
When participating in the sharing market, the optimization problem of each prosumer $i\in \mathcal{I}$ becomes
\bsq
\label{eq:SMK}
\bq
\mathop {\min }\limits_{{p_i},{b_i}} && \Pi_i:={c_i}p_i^2 + {d_i}{p_i} + (a{\lambda _c(b)} + {b_i}){\lambda _c(b)}\label{SMK.1}\\
{\rm{s}}.{\rm{t}}.&& {p_i} + a{\lambda _c(b)} + {b_i} = {D_i} \label{SMK.2}
\\
&&\sum\nolimits_i {(a{\lambda _c(b)} + {b_i})}  = Na{\lambda _c} + \sum\nolimits_i {{b_i} = } 0  \label{MCC}
\eq
\esq
Here, $\Pi_i(p_i,b_i,b_{-i})$ is the cost function, which can be divided into two parts: the disutility in terms of money $f_i(p_i):=c_ip_i^2+d_ip_i$ and the sharing cost $s_i(b):=(a\lambda_c+b_i)\lambda_c$. \eqref{SMK.2} represents the energy balancing. \eqref{MCC} is the sharing market clearing condition, which appears in every prosumer's problem. Due to the common constraint \eqref{MCC},  problem \eqref{eq:SMK} constitutes a generalized Nash game (GNG), where players' payoffs and strategy sets depend on each other.

In summary, the sharing game consist of the following elements: 1) the set of prosumers $\mathcal{I}=\{1,2,...,I\}$; 2) action sets $X_i(p_{-i},b_{-i})$\footnote{The subscribe $-i$ means all players in $\mathcal{I}$ except $i$ },$\forall i$, and strategy space  $X = \prod_i X_i$; 3) cost functions $\Gamma_i(p_i, b_i,p_{-i},b_{-i}), \forall i$. For simplicity, we  use $\mathcal{G}=\{\mathcal{I},X,\Gamma\}$ 
to denote the sharing game in an abstract form. 
\vspace{-1em}
\section{Properties of the Sharing Game}
\subsection{Existence and Uniqueness of Equilibrium}
In this subsection, we  show that the GNG model of the energy sharing problem can be reduced into a standard Nash game, and we prove the existence and uniqueness of its equilibrium.  

Denote by $b_j$  the bids of other prosumer $j$ $(j \ne i)$. From \eqref{MCC}, we have
\bq
\label{eq:lambda}
\lambda_c(b)=-\frac{b_i}{Na}-\frac{\sum \nolimits_{j \ne i} b_j}{Na}
\eq
Substituting  into \eqref{SMK.2} yields
\bq
\label{eq:p_i}
p_i= D_i-\frac{N-1}{N}b_i+\frac{\sum \nolimits_{j \ne i} b_j}{N}
\eq
Using $b_i$ to represent $p_i$ and $\lambda_c(b)$, the GNG \eqref{eq:SMK} degenerates into an equivalent standard Nash game \eqref{eq:Nash game}.
\bq
\label{eq:Nash game}
\mathop {\min }\limits_{{b_i}} && c_i \left(D_i-\frac{N-1}{N}b_i+\frac{\sum \nolimits_{j \ne i} b_j}{N}\right)^2 \nonumber\\
&&+ d_i \left(D_i-\frac{N-1}{N}b_i+\frac{\sum \nolimits_{j \ne i} b_j}{N}\right) \nonumber \\
&& + \left(-\frac{b_i}{N}-\frac{\sum \nolimits_{j \ne i}b_j}{N}+b_i\right)\left(-\frac{b_i}{Na}-\frac{\sum \nolimits_{j \ne i} b_j}{Na}\right)
\eq

Direct computation shows that, the second derivative of the objective function is $2\left[c_i\left(\frac{N-1}{N}\right)^2-\frac{N-1}{N^2a}\right]>0$, implying each prosumer solves a strictly convex optimization.

\begin{definition}(Nash Equilibrium)
	A strategy profile $ (p^*,b^*) \in X$ is a Nash Equilibrium (NE) of the sharing game $\mathcal{G}=\{\mathcal{I},X,\Gamma\}$\footnote{Given a collection of $x_i$ for $i$ in a certain set $A$, $x$ denotes the vector $x := (x_i; i\in A)$ of a proper dimension with $x_i$ as its components}, if $\forall i \in \mathcal{I}$ 
	\bq
	\Gamma_i(p_i^*,b_i^*,p_{-i}^*,b_{-i}^*) \le \Gamma_i(p_i,b_i, p_{-i}^*,b_{-i}^*), \forall (p_i,b_i) \in X_i(p_{-i}^*,b_{-i}^*) \nonumber
	\eq
\end{definition}

Specially, in the sharing game \eqref{eq:SMK}, the action set and the objective function of prosumer $i$ are independent of $p_{-i}$, so $X_i(p_{-i},b_{-i})$ reduces to $X_i(b_{-i})$, and $\Gamma_i(p_i,b_i,p_{-i},b_{-i})$ reduces to $\Pi_i(p_i,b_i,b_{-i})$. 
Given $p$, define $\tilde \lambda (p):=\frac{1}{N} \sum \nolimits_i (2c_ip_i+d_i)$ and $\tilde b_i(p):=D_i -p_i -a\tilde \lambda(p)$. We have the following proposition.
\begin{proposition}
	\label{Thm:prop-2}
	There exists a unique NE for the sharing game \eqref{eq:SMK}. Moreover, a strategy profile $(p^*,b^*)$ is the unique NE if and only if, $ \forall i \in \mathcal{I}$, $p^*_i$ is the unique solution of:
	\bsq
	\label{eq:Central}
	\begin{align} 
	\min_{p_i, \forall i}~~ & \sum\limits_i {\left({c_i} - \frac{1}{{2(N-1)a}}\right)p_i^2 + \left({d_i} + \frac{{{D_i}}}{{(N-1)a}}\right){p_i}}  \\
	\mbox{s.t.}~~ & \sum\nolimits_i {{p_i}}  = \sum\nolimits_i {{D_i}} :\xi 
	\end{align}
	\esq
	and $b_i^*=\tilde b_i(p^*)$.
\end{proposition}

The proof of Proposition \ref{Thm:prop-2}  can be found in Appendix A.
Proposition \ref{Thm:prop-2} is fundamental since it ensures that the proposed sharing game is well defined. Furthermore, it implies that the NE computation can be greatly simplified into solving  a  simpler optimization problem  \eqref{eq:Central}, which is  strictly convex. 
\vspace{-1em}
\subsection{Incentives for Prosumers' Participation}
The next proposition shows that all prosumers are incentivized to share by comparing the costs of the individual decision-making and the sharing game \eqref{eq:SMK} at equilibrium.   

Let $\Pi_i(p^*_i,b^*)$ be the cost of prosumer $i$ at the NE of sharing game $\mathcal{G}=\{\mathcal{I}, X,\Pi\}$ defined by \eqref{eq:SMK}, and $f_i(D_i)$ the cost of prosumer $i$ with his individual optimal decision.
\begin{proposition}
	\label{Thm:prop-3} We have
	\bq
	\label{eq:pareto}
	\Pi_i (p^*_i,b^*) \le f_i(D_i), \forall i\in \mathcal{I} 
	\eq
	moreover, \eqref{eq:pareto} holds with strictly inequality for at least one $i$ unless the unique optimal solution of \eqref{eq:Central} is $p_i^*=D_i, \forall i$.
\end{proposition}

The proof of Proposition \ref{Thm:prop-3} can be found in Appendix B. It says that with the proposed sharing mechanism, a Pareto improvement can be achieved for all prosumers, since the cost of each prosumer is no worse than making decisions individually.  Hence, the sharing mechanism provides positive incentives for prosumers to participate in the sharing market. It is worth noting that, a Pareto improvement is achieved for all prosumers via energy sharing, but it is not necessarily a Pareto optimum. In game theory, every player takes strategic action under rationality, resulting in a Nash equilibrium, which is usually not Pareto optimal \cite{nash1953two}. This does not contradict to our conclusion here.

\subsection{Sharing Price and Prosumers' Behavior}

In this subsection, we clarify the relationship between the sharing price $\lambda_c$ that clears the market and prosumer's marginal disutility, as well as the resulting prosumers' behaviors.

\begin{definition} (Marginal Disutility)
	The marginal disutility of prosumer $i$, denoted by $md_i$, is defined as
	\bq
	md_{i}(p_i):=\frac{\partial f_i(p_i)}{\partial p_i}=2c_ip_i+d_i \nonumber
	\eq
\end{definition}

\begin{proposition}
	\label{Thm:prop-1}
	Assume $(p^*,b^*)$ is the NE of the sharing game \eqref{eq:SMK}. Then
	\begin{enumerate}
		\item the  sharing price at equilibrium is given by
		\bq
		\lambda_c^* =\frac{1}{N} \sum \limits_i md_{i}(p_i^*); \nonumber
		\eq
		\item $md_{i}(p_i^*)>\lambda_c(b^*)$ if and only if $q_i(b^*)>0$.
	\end{enumerate}
\end{proposition}

The proof of Proposition \ref{Thm:prop-1} can be found in Appendix C.
Proposition \ref{Thm:prop-1} says that, the clearing price at the NE is simply the average marginal disutility of all prosumers participating in the sharing market. Moreover, the prosumer whose marginal disutility is larger than the average (which equals $\lambda_c$) has $q_i(b^*)>0$ and hence will buy energy, while whose marginal cost is lower than the average has $q_i(b^*)<0$ and hence will sell energy. Under the proposed sharing mechanism, a prosumer with higher/lower marginal disutility adjusts less/more and purchases/sells in the sharing market.


\subsection{Social Efficiency}
To investigate the social efficiency of the proposed sharing mechanism, consider the social planner's problem:
\bsq
\label{eq:AGG}
\bq
\mathop {\min }\limits_{{p_{i,\forall i \in \mathcal{I}}}} && \sum\nolimits_i {({c_i}p_i^2 + {d_i}{p_i})} \label{AGG.1}\\
{\rm{s}}{\rm{.t}}{\rm{.}} && \sum\nolimits_i {{p_i}}  = \sum\nolimits_i {{D_i}} \label{AGG.2}
\eq
\esq 

\begin{definition} (Socially Optimal)
	$\bar p$ is socially optimal if $\bar p$ is the unique optimal solution of \eqref{eq:AGG}.
\end{definition}

Optimal solution of \eqref{eq:AGG} is different from that under individual decision-making, except for the case in which $p_i=D_i$ $\forall i\in \mathcal{I}$ happens to be the optimal solution to problem \eqref{eq:AGG}. The difference in their optimal values interprets the loss of social welfare. Next we reveal that the proposed energy sharing mechanism can effectively reduce the loss of social welfare.

Comparison of \eqref{eq:Central} in Proposition \ref{Thm:prop-2} with the social planner's problem \eqref{eq:AGG} suggests that the NE of the sharing problem \eqref{eq:SMK} 
might approach social optimality as $N\rightarrow \infty$.  The next result states this intuition precisely.

\begin{proposition}
	\label{Thm:prop-4}
	Let $(p^*(N),b^*(N))$ be the unique NE of \eqref{eq:SMK}, and $\bar p(N)$ be the socially optimal solution of \eqref{eq:AGG}. Then, we have
	$$\sum \nolimits_{i \in \mathcal{I}} f_i(p_i^*(N)) \ge \sum \nolimits_{i \in \mathcal{I}} f_i(\bar p_i(N))$$
	and the average cost difference
	\bq
	\mathop {\lim }\limits_{N \to \infty } \frac{1}{N}\left[ {\sum\nolimits_{i \in \mathcal{I}} {{f_i}(p_i^*(N)) - \sum\nolimits_{i \in \mathcal{I}} {{f_i}({{\bar p}_i}(N))} } } \right] = 0 \nonumber
	\eq
\end{proposition}

The  proof of Proposition \ref{Thm:prop-4} can be found in Appendix D. The social planner's problem can be regarded as a utility maximization approach exerted by a non-profit utility. It minimizes the disutility of all prosumers.
Proposition \ref{Thm:prop-4} 
says that the proposed sharing mechanism asymptotically converges to the social optimum when there is a large enough number of prosumers involved. It is worth noting that Proposition \ref{Thm:prop-4} is essentially a sufficient condition. Furthermore, as shown in Proposition \ref{Thm:prop-2} that a unique Nash equilibrium exists in our sharing game, it can be proved that the price-of-anarchy (PoA)
\bq
{\rm{PoA}} =  \frac{\sum_{i \in \mathcal{I}} f_i(p_i^*(N))}{\sum_{i \in \mathcal{I}} f_i(\bar{p}_i (N))} \le 1+\frac{\beta}{N}
\eq
where $\beta$ is a positive number. The obtained PoA is in the same order of magnitude as given in \cite{johari2011parameterized}, which is $1+1/(N-2)$.

Similarly, the impact of price sensitivity can be analyzed by the following proposition.
\begin{proposition}
	\label{Thm:prop-4-2}
	Let $(p^*(a),b^*(a))$ be the unique NE of \eqref{eq:SMK}, with price sensitivity equals to $a<0$. Then, we have 
	$\sum \nolimits_{i \in \mathcal{I}} f_i(p_i^*(a))$
	is decreasing in $|a|$.
\end{proposition}

The proof of proposition \ref{Thm:prop-4-2} can be found in Appendix E. It  reveals that when the prosumers are more sensitive to the change of price, the total social cost under sharing decreases and becomes closer to the social optimal cost. 
It is worthy nothing that, because $\sum_{i\in \mathcal{I}} (a \lambda_c+b_i)\lambda_c=0$ holds, the group of prosumers are budget self-balancing, which is a main superiority compared with the VCG mechanism. 

\textbf{Remark 2:} Though currently there is no mature sharing platform for energy sharing in practice yet, it could potentially be provided in the following ways.

1) \emph{Offered by the utility operator}. It is proved in Proposition 2 that the cost of each prosumer will never be worse than when the prosumer does not share so that every prosumer has an incentive to participate in sharing. At least one prosumer has a positive cost reduction except for the extreme case where the individual decisions happen to be the same as the social optimum. It means that $\Delta f := \sum_{i \in \mathcal{I}} f_i(D_i)- \sum_{i \in \mathcal{I}} f_i(p_i^*)>0$. In this context, part of this cost reduction $\alpha \Delta f (0<\alpha<1)$ can be taken as a payment to the utility operator and each prosumer $i$ pays $\alpha[f_i(D_i)-f_i(p_i^*)]$.

2) \emph{Offered by a third-party company}. In sharing economy, sometimes a third-party company is willing to provide such a platform free of charge. Instead of earning money directly from the prosumers, they make money by other ways such as putting advertisements on the platform.

3) \emph{Offered by a non-profit entity}. In some countries, the power grid is operated by a state-owned company, which aims to enhance social welfare, other than its profit. As shown in the paper, a Pareto improvement is achieved via energy sharing. It seems reasonable that this non-profit entity may have the incentive to provide such a platform.
\section{Impacts of Competition}
Above analysis assumes a perfectly monopolistic competition market, where each prosumer owns only one resource. Next we analyze a more complicated situation, in which  prosumers could own multiple resources. 

Assume that there are $I$ multi-resource prosumers (MRPs) indexed by $i \in \mathcal{I}=\{1,2,...,I\}$. Each prosumer $i$ owns $K_i$ kinds of resources labeled by $k \in \mathcal{K}_i = \{1,2,...,K_i\}$. The corresponding 
energy productions are $p_i^1, p_i^2, \cdots, p_i^{K_i}$.  However, we assume there are still $N$  resources in total, which means $\sum \nolimits_{i\in \mathcal{I}} K_i =N$. 
The multi-resource prosumer (MRP) $i\in \mathcal{I}$ can either carry out the demand response command individually by solving the following problem
\bsq
\label{eq:Ind-Ex}
\bq
\mathop {\min }\limits_{{p_i^k, \forall k  \in \mathcal{K}_i}} && \sum \nolimits_{k} \left[{c_i^k}\left(p_i^k\right)^2 + {d_i^k}\left(p_i^k\right)\right]
\label{Ind-Ex.1}\\ 
{\rm{s}}{\rm{.t}}{\rm{.}}  && \sum \nolimits_{k} {p_i^k} = {D_i}\label{Ind-Ex.2}
\eq     
\esq
or take part in the sharing market by solving
\bsq
\label{eq:SMK-Ex}
\bq
\mathop {\min }\limits_{{p_i^k,\forall k \in \mathcal{K}_i},{b_i}} && \sum\nolimits_k \left[{c_i^k}\left(p_i^k\right)^2 + {d_i^k}{p_i^k}\right] + (a{\lambda _c} + {b_i}){\lambda _c}\label{SMK-Ex.1}\\
{\rm{s}}.{\rm{t}}.&& \sum\nolimits_k {p_i^k} + a{\lambda _c} + {b_i} = {D_i} \label{SMK-Ex.2} \\
&& \sum \nolimits_i (a{\lambda _c}) + \sum\nolimits_i {{b_i}}  = 0 \label{SMK-Ex.3}
\eq
\esq
Following the similar process as in Section III, we can easily prove that the proposed sharing mechanism can benefit all MRPs and the equilibrium sharing price reflects the average marginal disutility.

\begin{definition} (Social Optimal for MRP) $\bar p$ is socially optimal (or the most efficient) if $\bar p$ solves
	\bsq
	\label{eq:AGG-Ex}
	\bq
	\mathop {\min }\limits_{{p^k_i, \forall i \in \mathcal{I}, k \in \mathcal{K}_i}} && \sum\nolimits_i \sum\nolimits_k {\left[{c_i^k}\left(p_i^k\right)^2 + {d_i^k}{p_i^k}\right]} \label{AGG-Ex.1}\\
	{\rm{s}}{\rm{.t}}{\rm{.}} && \sum\nolimits_i  \sum\nolimits_k {{p_i^k}}  = \sum\nolimits_i {{D_i}} \label{AGG-Ex.2}
	\eq
	\esq  
	
\end{definition}

If $I=1$, the  sharing problem for MRP \eqref{eq:SMK-Ex} becomes the social optimal problem \eqref{eq:AGG-Ex}. If $I=N$, the sharing problem for MRP \eqref{eq:SMK-Ex} degenerates to the sharing problem in the previous sections. If $1<I<N$, the following propositions hold. 

\begin{proposition}
	\label{Thm:prop-7}
	The marginal disutilities of the resources for a prosumer are equal, which means\footnote{$p_i$ denotes the vector $p_i:=(p_i^k, \forall k)$}
	\bq
	md_{i}(p_i ):=md_{i}^1\left(p_i^1\right)=\cdots=md_{i}^{K_i}\left(p_i^{K_i}\right),\quad \forall i\in \mathcal{I}
	\eq
	
\end{proposition}
Proposition \ref{Thm:prop-7} can be directly deduced from the KKT condition. So we omit the proof here.

Given $p:=(p_i,\forall i)$, we define $\tilde \lambda (p):=\frac{1}{I} \sum \nolimits_i md_{i}(p_i)$ and $\tilde b_i(p):=D_i -\sum \nolimits_k p_i^k -a\tilde \lambda(p)$. Then we have the following proposition.
\begin{proposition}
	\label{Thm:prop-6}
	There exists a unique NE for the sharing problem with MRP \eqref{eq:SMK-Ex}. Moreover, a strategy profile $(p^*,b^*)$ is the unique NE if and only if, $\forall i\in \mathcal{I}$, $p^*_i$ is the unique solution of \eqref{eq:Central-Ex}.
	\bsq
	\label{eq:Central-Ex}
	\begin{align}
	\mathop {\min }\limits_{{p^k_{i}, \forall i, k \in \mathcal{K}_i}}~~ & \sum\nolimits_i \sum\nolimits_k \left[\left({c_i^k}-\frac{1}{2(I-1)a}\right)\left(p_i^k\right)^2 \right.\nonumber \\
	& + \left.\left({d_i^k}+\frac{D_i}{(I-1)a}\right){p_i^k}\right ]-\frac{\sum\limits_i \sum\limits_{j>k \in \mathcal{K}_i} p_i^k p_i^j}{(I-1)a}  \label{Central-Ex.1} \\
	\mbox{s.t.}~~ & \sum\nolimits_i \sum\nolimits_k {{p_i^k}}  = \sum\nolimits_i {{D_i}} :\xi^\prime \label{Central-Ex.2}
	\end{align}
	\esq
	and $b_i^*=\tilde b_i(p^*)$.
\end{proposition}

The proof of Proposition \ref{Thm:prop-6}  can be found in Appendix F. Proposition \ref{Thm:prop-6} extends the result of existence and uniqueness of the NE in the sharing game from the single-resource case to the multi-resource one, and again provides an effective way to simplify the computation of NE. 

Then we analyze the change in efficiency based on model \eqref{eq:Central-Ex}. Generally speaking, as $I$ varies from $1$ to $N$, the change in the total social optimal cost may not be monotonous. So we only consider a special case in which all $c^k_i=c, \forall k, \forall i$ and $K_i=K_I, \forall i$ and gives the following proposition.

Let $(I,K_I,D)$ denotes a scenario that there are $I$ prosumers, each has $K_I$ resources and the required load adjustment for prosumer $i$ is $D_i$. Then, the scenario $(I^{'},K_{I^{'}},D^{'})$ is an \textbf{\emph{equal partition}} of $(I,K_I,D)$ when there exists an $Z \in \mathbb{Z^{+}}$, such that $I^{'}=ZI$ and $K_I=ZK_{I^{'}}$, the resources one prosumer possesses and required load adjustment is distributed equally to $Z$ prosumers and satisfies, $\forall i, \forall z_1, z_2 \in \{1,...,Z\} $
\bq
\label{condition}
({D_{Z(i-1)+z_1}^{'}} - \sum\limits_{k=1}^{K_{I^{'}}} {p_{_{Z(i-1)+z_1}}^{k*}} )({D_{Z(i-1)+z_2}^{'}} - \sum\limits_{k=1}^{K_{I^{'}}} {p_{_{Z(i-1)+z_2}}^{k*}} )  \ge  0 \nonumber\\
\sum \nolimits_{z=1}^Z D_{{Z(i-1)+z}}^{'} =  D_i 
\eq
where $p^*$ corresponds to the NE under scenario $(I,K_I,D)$.
\begin{definition}(Variance of marginal disutility) The variance of marginal utilities $md_i^{k*}, \forall k\in \mathcal{K}_I, \forall i\in \mathcal{I}$ is defined as
	$$Var(md_i^{k*}(I),I):=\frac{1}{N} \sum \limits_{i=1}^{I} \sum \limits_{k=1}^{K_I} (md_{i}^{k*}-\frac{1}{N} \sum \limits_{i=1}^{I} \sum \limits_{k=1}^{K_I} md_{i}^{k*})^2$$
	\label{def-5}
\end{definition}

\begin{proposition}
	\label{Thm:prop-8}
	Suppose $c_i^k=c, \forall i\in \mathcal{I}, \forall k\in \mathcal{K}_i$ and $(p^*(I),b^*(I))$ is the unique NE of the sharing problem for MRP \eqref{eq:SMK-Ex}  with $I>1$. For any $I$ prosumers with the same number of resources, i.e. $K_i=K_I, \forall i\in\mathcal{I}$, there always exists an equal partition of $(I,K_I,D)$ to $(I^{'},K_{I^{'}},D^{'})$, such that
	$$\sum \nolimits_{i=1}^{I^{'}} \sum \nolimits_{k=1}^{K_{I^{'}}} f_{ik} (p_i^{k*}(I^{'})) \le \sum \nolimits_{i=1}^{I} \sum \nolimits_{k=1}^{K_I} f_{ik} (p_i^{k*}(I))$$
	Moreover, if $-2ac \le N$, then
	$$ Var(md_i^{k*}(I^{'}), I^{'}) \le Var(md_i^{k*}(I), I) $$
\end{proposition}
The proof can be found in Appendix G. It shows that the system under $I=1$ is the most efficient; otherwise, introducing competition by spreading resources benefits social welfare. Proposition \ref{Thm:prop-8} only considers a very special case. In Section IV, we provide empirical results of numerical experiments to further confirm this property.  

\section{Possible Extensions}
This paper considers a simplified situation (power balance) to study the strategic behavior of individual prosumers in a sharing scheme.  We remark on
some possible extensions.

1) \textbf{Network constraints.} Our model is restricted to a community with a few buildings or microgrids that accounts for only a small fraction of the total demand at a node of a city-sized distribution network. In such a situation, it is reasonable to neglect network constraints as our system is smaller than a node, as done in \cite{wang2016cooperative,choi2016advanced,mondal2015distributed}.
To incorporate prosumers with larger capacities that can take part in a distribution market, one can replace condition \eqref{eq:clearcondition} with an optimal power flow based market clearing problem. This is among our undergoing works.

2) \textbf{Uncertainty.} If a prosumer has a random generation device, assume its output is $\tilde{p}_i^r=p_i^{r0}+\Delta p_i^r$, where $p_i^{r0}$ is the forecast output and $\Delta p_i^r$ the uncertain deviation. In the day-ahead stage, $E_i^0$ is determined according to the forecast value $p_i^{r0}$ and we have $p_i^0+E_i^0=\hat{D}_i^0$ (where $\hat{D}_i^0:=D_i^0-p_i^{r0}$). Then in the real-time stage, after the exact output of random devices and the load reduction requirement $D_i$ are known, prosumer $i$ adjusts the output of his controllable resource by $p_i$ to maintain energy balance, which means $p_i^0+p_i+(E_i^0-D_i)=\hat{D}_i^0-\Delta p_i^r$. So we have $p_i=D_i-\Delta p_i^r$. Since $\Delta p_i^r$ is uncertain, the stochastic optimization approach can be applied to solve this problem by replacing the objective with the expectation of his cost. By simply approximating the expectation by a weighted sum of objective functions under different scenarios, a similar procedure to that in this paper can be used check that the main conclusions still hold.

3) \textbf{Heterogeneous participants.} To accommodate heterogeneous participants, \eqref{eq:clearcondition1}-\eqref{eq:clearcondition} can be regarded as set rules in the sharing market and some constraints can be added to depict the participant's feature. For instance, if some of the prosumers are aggregated, it is equivalent to one aggregated prosumer owning various resources as in model \eqref{eq:SMK-Ex}. If the participant has storage facilities, related constraints \cite{Sharing-analytical1} can be included. If the participant has limited adjustable capacity, bound constraint $\underline{p}_i \le p_i \le \overline{p}_i$ can be added. If, on the other hand, prosumers have different but fixed $a_i$ and still bid on only $b_i$, then the sharing game generalizes directly with the sharing price $\lambda_c = -\sum_i b_i / \sum_i a_i$.
Moreover, it can be shown that Propositions \ref{Thm:prop-2}--\ref{Thm:prop-1} hold with straightforward modifications.
In Section VI.B we present simulation results that suggest that the Nash equilibrium may still converge to the social optimal even with heterogeneous $a_i$
(cf Proposition \ref{Thm:prop-4}).

4) \textbf{Bounded rationality.} In our model, each prosumer is connected to the platform via a smart meter and the smart meter bids automatically. In this context, it is reasonable to assume that the smart meters' algorithm is designed with full rationality as long as the information required in computation is available. Still, it might be more practical in some cases where bounded rationality should be considered \cite{saad2016toward} and prospect theory can be adopted to characterize this phenomenon. Prospect theory extends the traditional expected utility theory by rewriting the objective function taking into account the following two aspects: i)\emph{Weighting Effect.} It has been found that most people overweight/underweight low/high probability income . To take this into account, instead of an objective probability $p$ in the stochastic model in 2) above, a subjective weighting function $W(p)$ can be incorporated. ii) \emph{Framing Effect.} People prefer stable incomes, but when facing losses, they may prefer taking risks to compensate for losses. The initial wealth will influence people's utility as well. To take this into account, a reference point is chosen; the coefficient $c_i,d_i$ in the model are not necessarily the objective cost parameters, and the prosumer can input his subjective parameters instead. 

\section{Illustrative Examples}

In this section, numerical experiments are presented to illustrate theoretical results. First, a simple case is used to clarify the basic setup. Then, the impacts of several factors are analyzed, including the number of prosumers, their price sensitivities as well as  the impact of competition.

\subsection{Benchmark Case}
\begin{figure}[!t]
	\centerline{\includegraphics[width=0.7\columnwidth]{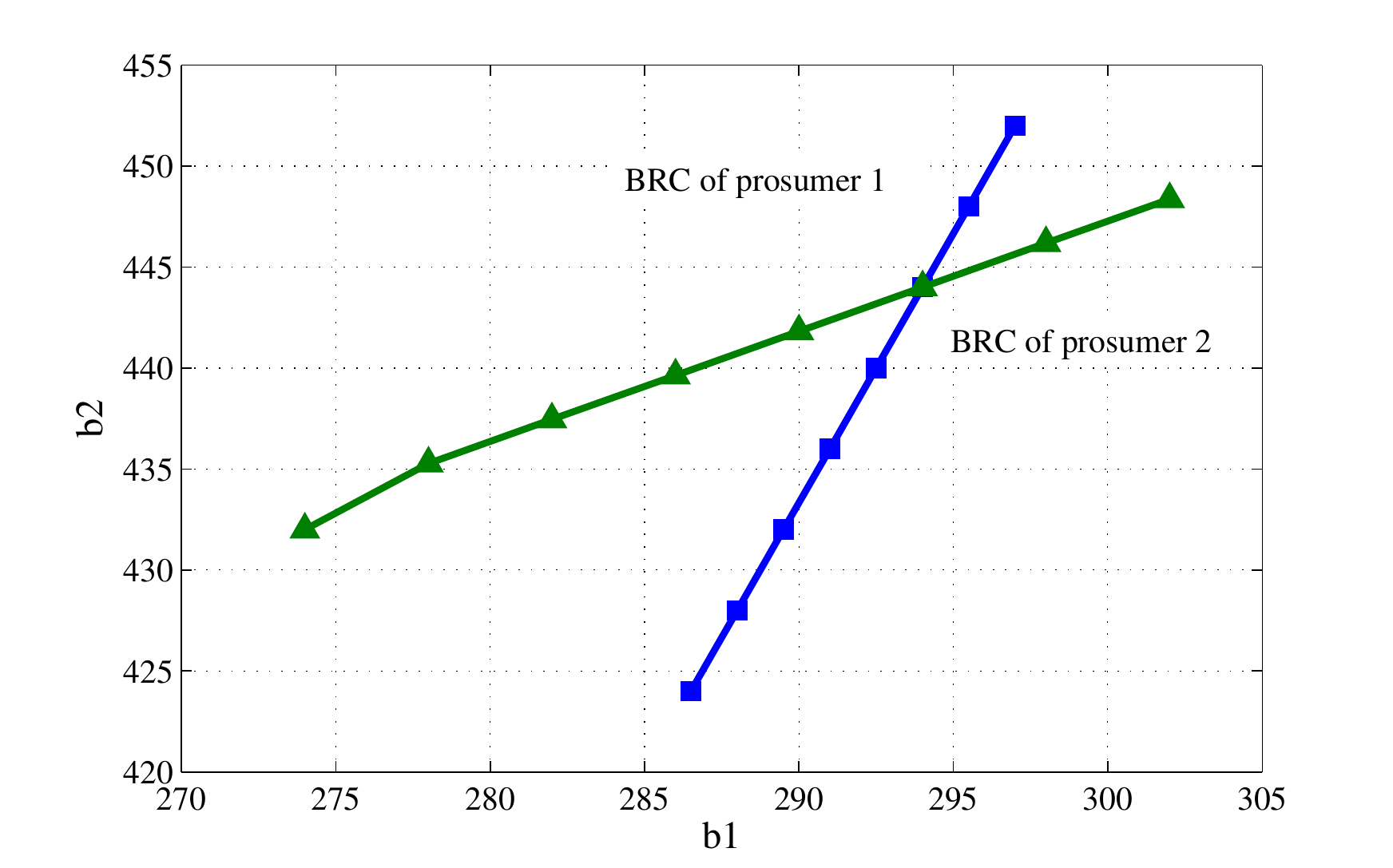}}
	\setlength{\abovecaptionskip}{0pt}
	\caption{Best response curves of two prosumers.}
	\setlength{\belowcaptionskip}{-5em}
	\vspace{-1.5em}
	\label{fig:BRC}
\end{figure}

The simplest scenario with two prosumers is taken as an illustrative example. $p_1$ and $p_2$ are the output adjustment of prosumer 1 and 2. We assume the price sensitivity $a=-200$, the cost coefficients $c_1=0.003 \$/kW^{2}$, $d_1=0.042 \$/kW$ and $c_2=0.006 \$/kW^{2}$, $d_2=0.072 \$/kW$. The required demand reduction are $D_1=100 kW$ and $D_2=200 kW$. The optimal output adjustments and the corresponding costs when making decisions individually (IDL), taking part in the sharing market (SMK) and the social optimal (SCO) are shown in Table \ref{tab:benchmark}. The best response curves of two prosumers are shown in Fig.\ref{fig:BRC}.
\begin{table}[htp]
	\tiny
	\vspace{-1em}
	\renewcommand{\arraystretch}{1.3}
	\renewcommand{\tabcolsep}{1em}
	\centering
	\caption{Optimal solution under IDL, SMK and SCO}
	\vspace{-1em}
	\label{tab:benchmark}
	\begin{tabular}{cccc}
		\hline 
		Scheme & IDL & SMK & SCO \\
		\hline
		Optimal output adjustment $p_1$ (kW) & 100 & 175 & 216.67 \\
		Optimal output adjustment $p_2$ (kW) & 200 & 125 & 83.33\\
		Cost of prosumer 1 (\$) & 72 & 27 & 231.83\\
		Cost of prosumer 2 (\$) & 384 & 322.13 & 101.67 \\
		Social total cost (\$) & 456 & 349.13 & 333.50 \\
		Relative cost difference & 36.73\% & 4.68\% & -- \\
		\hline
		\vspace{-1em}
	\end{tabular}
\end{table}

From Table \ref{tab:benchmark}, we can find that when the prosumers take part in sharing, their individual costs all decrease (prosumer 1 from \$72 to \$27, and prosumer 2 from \$384 to \$322.13), so does the social total cost, confirming Proposition \ref{Thm:prop-3}. The relative social cost difference \footnote{Relative social cost difference(IDL,SCO)=$\frac{\sum_{i \in \mathcal{I}} f_i(D_i)-\sum_{i \in \mathcal{I}}f_i(\bar p_i)}{\sum_{i \in \mathcal{I}}f_i(\bar p_i)}$ , Relative social cost difference(SMK,SCO)=$\frac{\sum_{i \in \mathcal{I}} f_i(p_i^*)-\sum_{i \in \mathcal{I}}f_i(\bar p_i)}{\sum_{i \in \mathcal{I}}f_i(\bar p_i)}$} between IDL and SCO is 36.73\% while between SMK and SCO is 4.68\%, showing that the sharing mechanism can greatly reduce the social total cost. The intersection of best response curves in Fig.\ref{fig:BRC} gives the sharing market equilibrium, which is $(b_1,b_2)=(294, 444)kW$ and the corresponding equilibrium output adjustment is $(p_1,p_2)=(175, 125)kW$, which is the same as the results in Table \ref{tab:benchmark} offered by the proposed equivalent model \eqref{eq:Central}, verifying Proposition \ref{Thm:prop-2}. Compared with the ``two-sided market'', our approach determines the role of a prosumer according to his willingness. Prosumer 1 becomes a seller because $b_1<\bar{b}$ while prosumer 2 as a buyer since $b_2>\bar{b}$. Compared with the ``single-sided market with set-price'', in our model, the platform only provides the means for prosumers to determine their own decisions and the sharing price.

\subsection{Impact of the Number of Prosumers}
\label{subsec:Nprosumers}
\begin{figure*}[!t]
	\centering
	\subfigure[Relative cost difference]{\includegraphics[width=0.65\columnwidth]{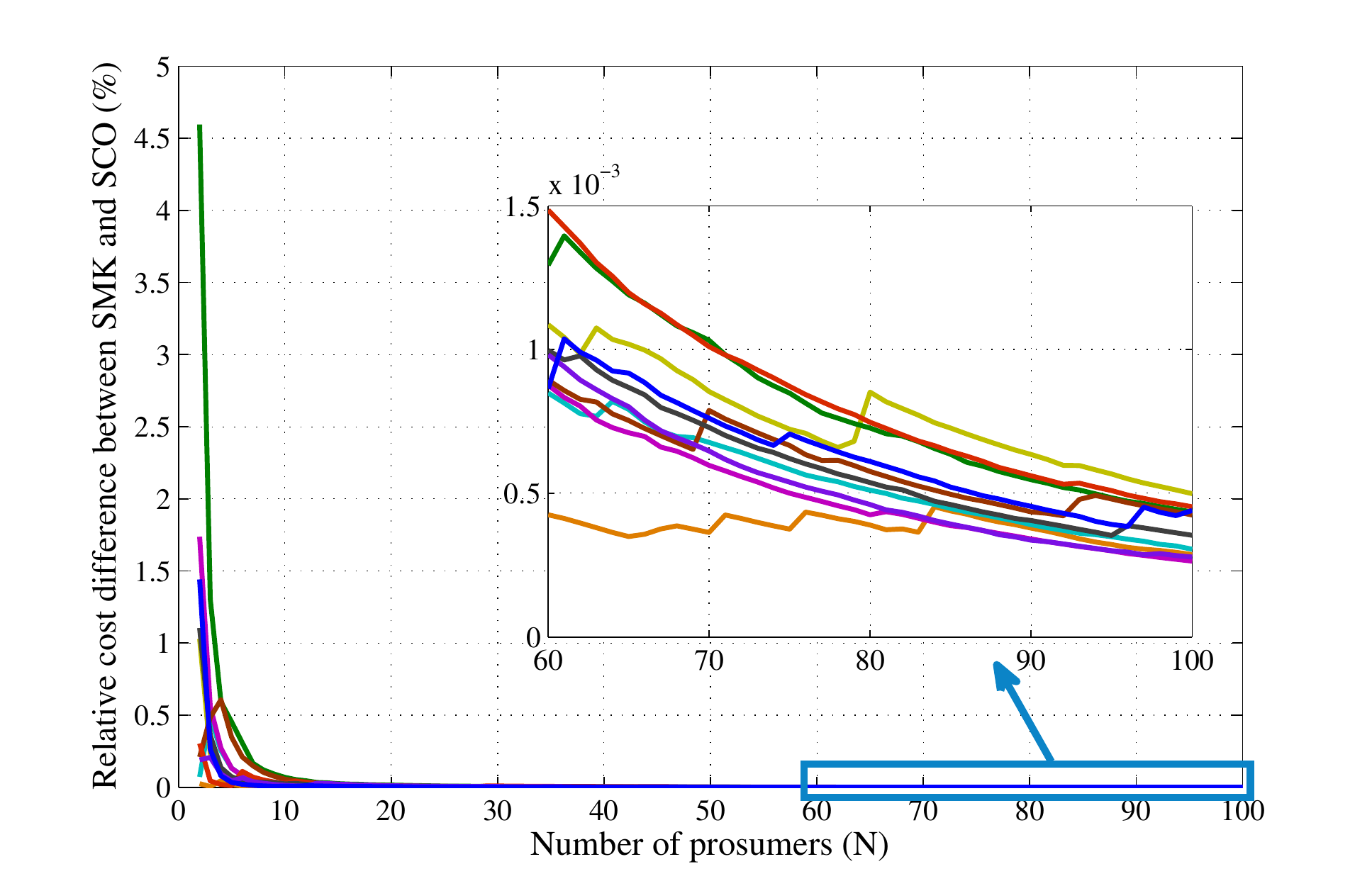}}
	\subfigure[Variance of marginal disutilities]{\includegraphics[width=0.65\columnwidth]{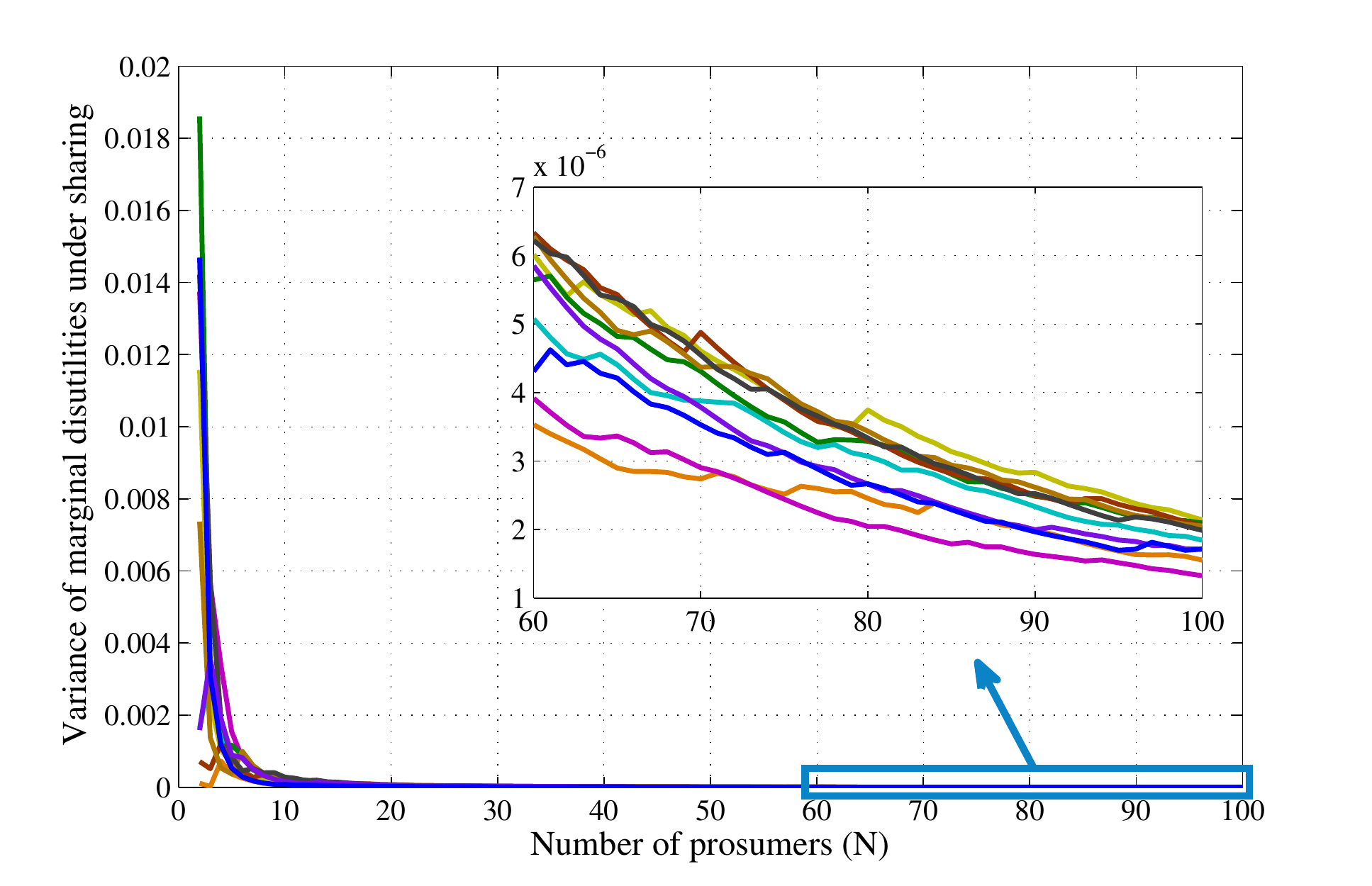}}
	\subfigure[Gap of average costs under different degree of heterogeneity]{\includegraphics[width=0.65\columnwidth]{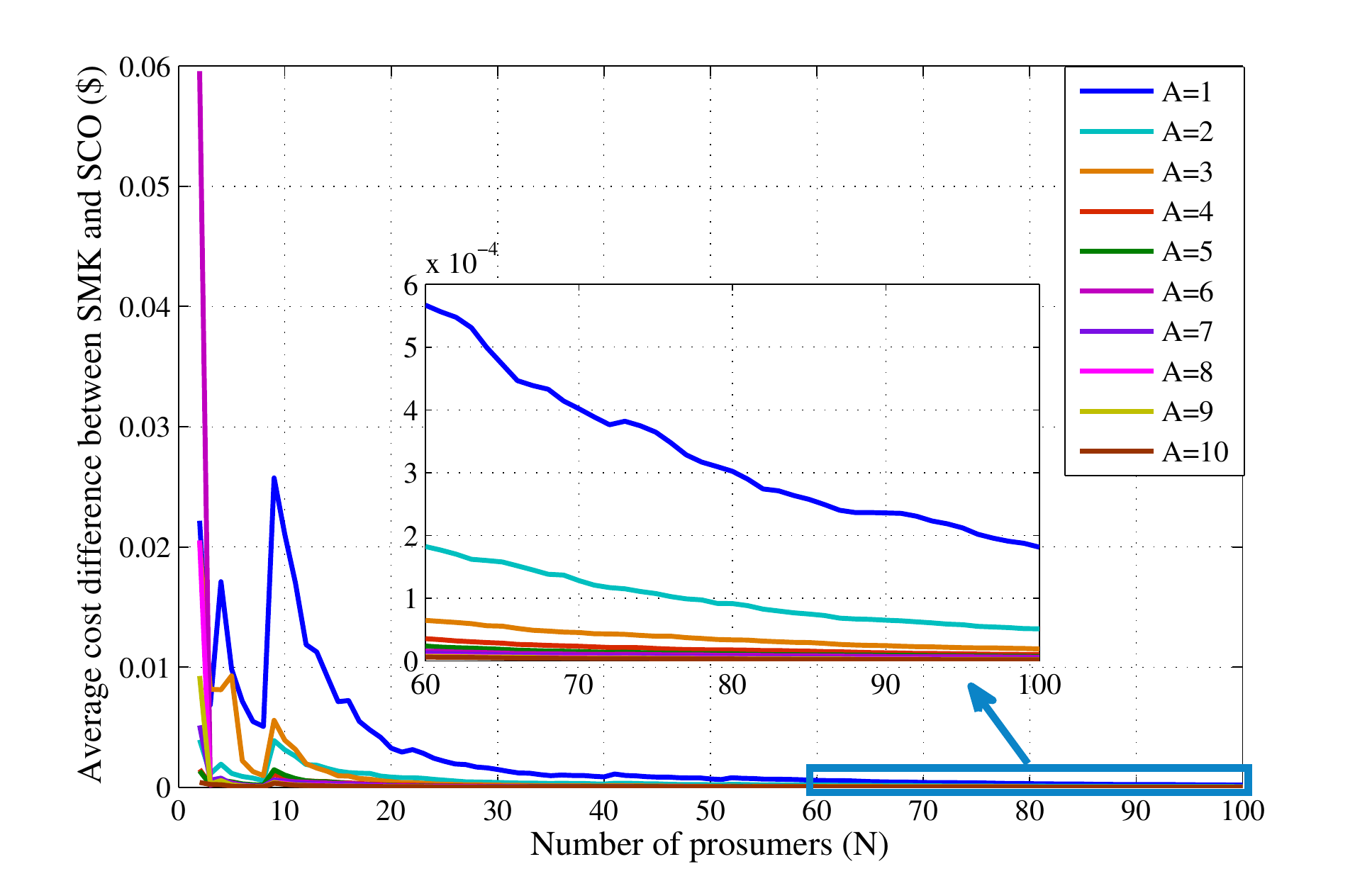}}
	\vspace{-1em}
	\caption{Results under different $N$.}
	\label{fig:different-N}
	\vspace{-1em}
\end{figure*}

We assume that $\underline c=0.001\$/kW^{-2}$, $\overline c=0.01\$/kW^{-2}$, $\underline d=0.02\$/kW$, $\overline d=0.12\$/kW$, $\underline D=0kW$, $\overline D=1000kW$ and change the number $N$ from 2 to 100. The parameters $c_i$, $d_i$, $D_i$ are randomly chosen within the upper and lower bounds and 10 scenarios are tested. For each of the 10 random scenario, the relative cost difference\footnote{The relative cost difference$:=[\sum_i f_i(p_i^*(N))-\sum_i f_i(\bar{p}_i(N))]/\sum_i f_i(\bar{p}_i(N))$} is plotted in Fig.\ref{fig:different-N}.(a) and the variance of marginal disutilities\footnote{The variance of marginal disutilities$:=\frac{1}{N}\sum_i(md_i-\sum_i md_i)^2$} in Fig. \ref{fig:different-N}.(b), both as functions of $N$. 

In Fig.\ref{fig:different-N}.(a), the relative cost differences are all positive, showing that the cost with sharing is always larger than the optimal social cost. When $N$ is larger than 60, the relative gap is less than $1.5\times 10^{-3}$ and the performance of sharing enhances with the number of participants. Moreover, the variance of marginal disutilities also drops sharply with increasing $N$, implying that the marginal disutilities of all prosumers become closer and all prosumers converge to the social optimum, as shown in Fig.\ref{fig:different-N}.(b). Actually, the relative gap is less than 5\% (the maximum appears at $N=2$) and shrinks sharply with the increase of $N$ in all scenarios, which indicates that the sharing can achieve a near-optimal solution even with a  small $N$.

We further run some simulations with heterogeneous $a_i$ to show that the key qualitative property, Proposition \ref{Thm:prop-4}, approximately holds. The parameters $a_i,\forall i$ are randomly chosen within the region $[Aa_0,0]$, where $a_0$ is the value of $a$ in the benchmark case. $A$ can be used to depict the degree of heterogeneity and we change it from 1 to 10. Under each $A$ and the chosen $a_i,\forall i$, the average cost difference\footnote{The average cost difference$:=\frac{1}{N}[\sum_i f_i(p_i^*(N))-\sum_i f_i(\bar p_i(N))]$} between the SMK and the SCO with different $N$ is plotted in Fig.\ref{fig:different-N}.(c). It can be found that the gaps always converge to zero, though the rate of decrease varies a little.

\begin{figure}[!t]
	\centerline{\includegraphics[width=0.7\columnwidth]{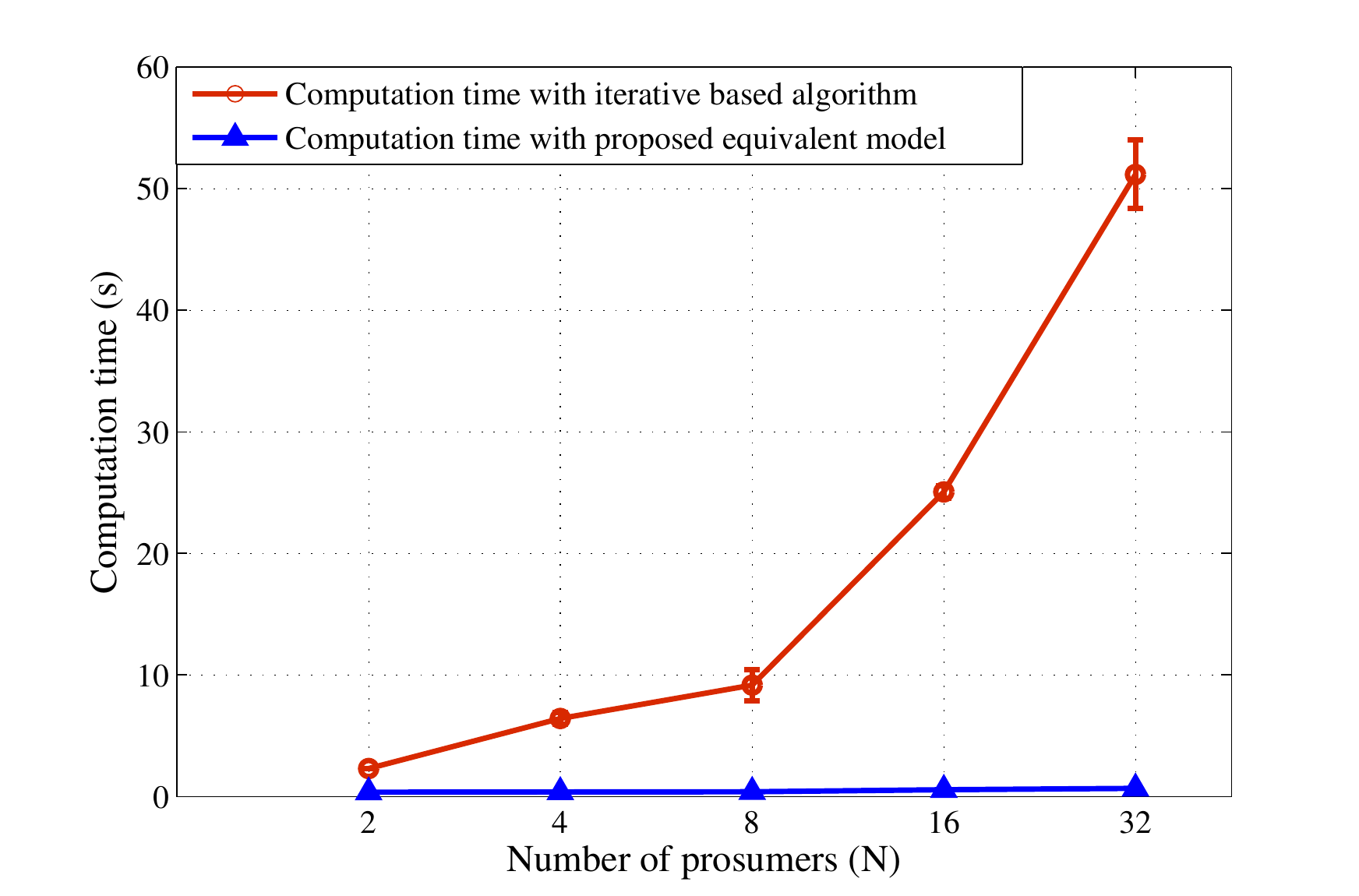}}
	\setlength{\abovecaptionskip}{0pt}
	\caption{Comparison of computation times with different $N$.}
	\setlength{\belowcaptionskip}{-5em}
	\vspace{-1em}
	\label{fig:computationaltime}
\end{figure}

The computation time using the proposed model \eqref{eq:Central} is compared with the traditional iterative based algorithm, which solves each prosumer's optimization problem sequentially with other prosumers' strategies fixed until convergence. The parameters $c_i$, $d_i$, $D_i,\forall i$ are randomly chosen within the upper and lower bounds and the scenarios with $N=2,4,8,16,32$ are tested. As shown in Fig.\ref{fig:computationaltime}, when the number of prosumers grows, the computation time using iterative based algorithm increases sharply, while the time required by the proposed algorithm changes little, showing the advantage of our method.

\subsection{Impact of Price Elasticity}
\begin{figure}[!t]
	\centerline{\includegraphics[width=0.7\columnwidth]{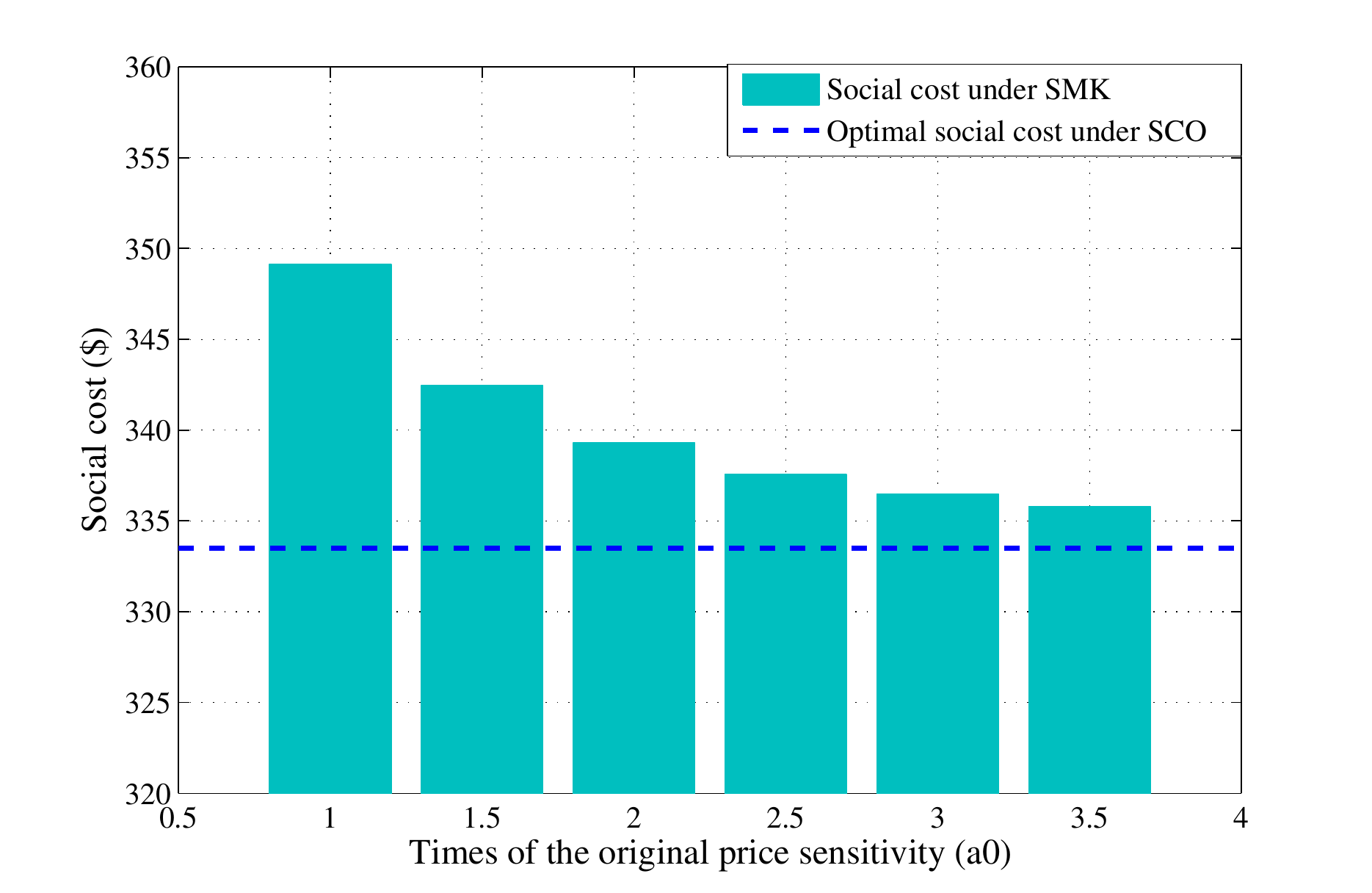}}
	\setlength{\abovecaptionskip}{0pt}
	\vspace{-0.5em}
	\caption{Social costs under different $a$.}
	\label{fig:price-elasticity}
	\vspace{-1.5em}
\end{figure}

When price sensitivity coefficient $a$ varies in $[1,3.5]$ times of its original value, the social costs under NE $\sum_i f_i(p_i^*)$ are shown in Fig. \ref{fig:price-elasticity}. The optimal social cost $\sum_i f_i(\bar{p}_i)$ is marked by a dash line. From the figure, when the absolute value of $a$ increases, the social cost under sharing is decreasing and gets closer to the optimal social cost. This is in accordance with Proposition \ref{Thm:prop-4-2}.
\subsection{Impact of Competition}
\begin{figure}[!t]
	\begin{minipage}[t]{0.5\linewidth}
		\centering
		\includegraphics[width=1.9in]{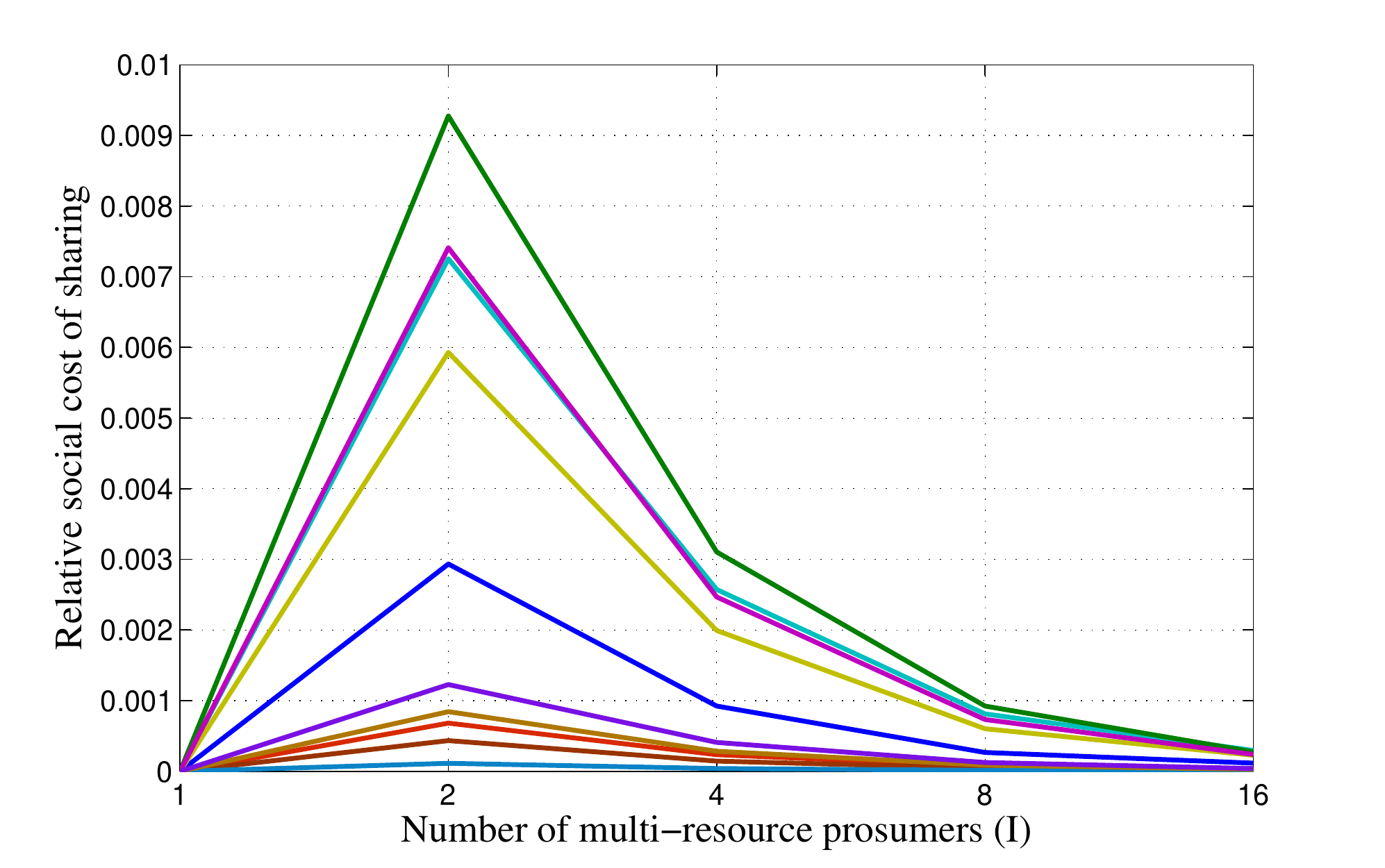}
		\centerline{\fontsize{8.5pt}{8.5pt}${\rm{(a)\; Same}} \;K_i$}
	\end{minipage}%
	\begin{minipage}[t]{0.5\linewidth}
		\centering
		\includegraphics[width=1.9in]{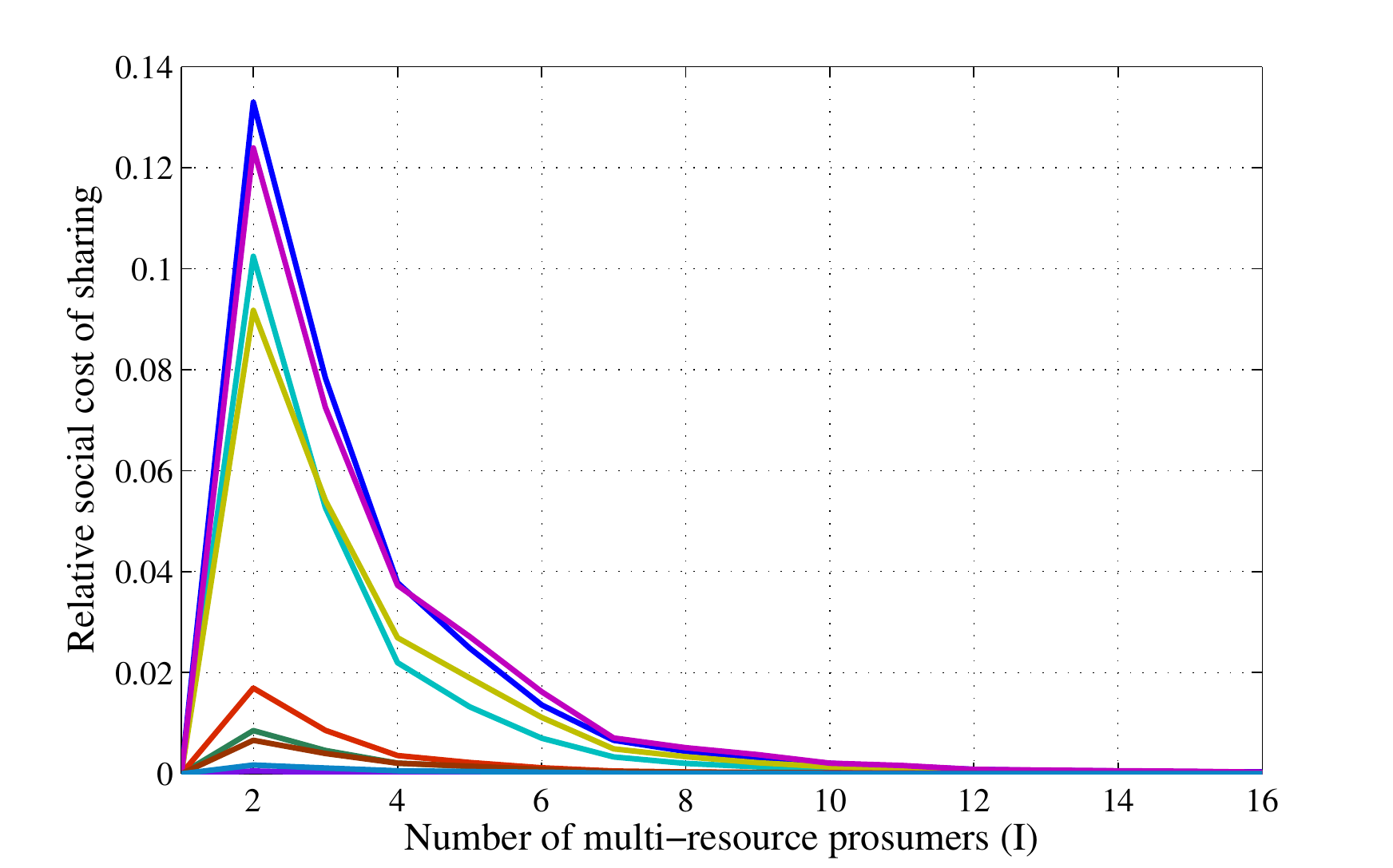}
		\centerline{\fontsize{8.5pt}{8.5pt}${\rm{(b)\; Different}} \;K_i$}
	\end{minipage}
	\caption{Relative cost of sharing  different numbers of MRPs with different $c_i^k$ and same/different $K_i$.}
	\label{fig:MRP}
\end{figure}

A special case, in which all $c^k_i=c, \forall k, \forall i$ and $K_i=K_I,\forall i$, is analyzed in Section IV, showing that introducing competition improves social welfare. However, the general case is difficult to prove. Here, we test cases with different $c_i^k$ but the same $K_i$; and cases with different $c_i^k$ and $K_i$. We assume that $\underline c=0.001\$/kW^{-2}$, $\overline c=0.01\$/kW^{-2}$, $\underline d=0.02\$/kW$, $\overline d=0.12\$/kW$, $\underline D=0kW$, $\overline D=1000kW$. The parameters $c_i^k$, $d_i^k$ and $D_i$ are randomly chosen within the ranges and 10 scenarios are tested. The resources are allocated according to the condition \eqref{condition} and the way in \cite{Appendix} (for same/different $K_i$). To eliminate the impact of scale, relative social cost, which equals to the ratio of costs with $I>1$ and $I=1$ minus 1,\footnote{Relative social cost$:=\frac{\sum\nolimits_{i=1}^{I} \sum \nolimits_{k=1}^{K_I} f_{ik}(p_i^{k*}(I))}{\sum\nolimits_{i=1}^{I^{'}=1} \sum \nolimits_{k=1}^{K_{I^{'}}} f_{ik}(p_i^{k*}(I^{'}=1))}-1$} is used for comparison and its change under different scenarios with different $c_i^k$ and same/different $K_i$ are given in Fig.\ref{fig:MRP}.(a) and Fig. \ref{fig:MRP}.(b), respectively.

It can be observed from Fig.\ref{fig:MRP}.(a) that the relative cost is the smallest when $I=1$, which means all the resources are owned by one prosumer and the social optimal is achieved. When $I \ge 2$, the relative cost decreases with $I$, demonstrating that competition improves economic efficiency as stated in Proposition \ref{Thm:prop-8}. This property extends to the case when $K_i$ are different as shown in Fig. \ref{fig:MRP}.(b).

\section{Conclusion}
By allowing prosumers to share energy with each other, the total cost of prosumers can be reduced, enhancing the social welfare. To promote energy sharing in smart grid, a simple mechanism is proposed, based on a generic supply-demand function. It allows prosumers to freely decide whether they will buy or sell based on their purchase desire. We have proved that a unique Nash equilibrium exists and the sharing price in the Nash equilibrium is equal to the average marginal disutility of the prosumers. We have proved that the proposed sharing scheme provides incentives for individual prosumers to participate and that the sharing system converges to the social optimal as the number of consumers increases. Case studies have been presented to illustrate the effectiveness of the proposed sharing mechanism. Future research directions include incorporating network constraints for larger scale systems, analyzing the behavior of this mechanism in the presence of uncertainty due to the integration of renewables and characterizing the impact of bounded rationality.

\ifCLASSOPTIONcaptionsoff
\newpage
\fi

\bibliographystyle{IEEEtran}
\bibliography{IEEEabrv,mybib}

\appendices

\makeatletter
\@addtoreset{equation}{section}
\@addtoreset{theorem}{section}
\makeatother
\renewcommand{\theequation}{A.\arabic{equation}}
\renewcommand{\thetheorem}{A.\arabic{theorem}}

\section{Proof of Proposition \ref{Thm:prop-2}}
\begin{proof}
	Given the other prosumers' bids $b_j, j \ne i$, prosumer $i$ solves problem \eqref{eq:SMK} to obtain his optimal strategy, which is equivalent to solving problem \eqref{eq:Nash game} with $p_i^*$ given by \eqref{eq:p_i}. Problem \eqref{eq:Nash game} is a strictly convex optimization problem. Thus, the KKT condition below is the necessary and sufficient condition for the optimal solution.
	\bq
	\label{eq:Sharing-KKT}
	2c_i(D_i-\frac{N-1}{N}b_i+\frac{\sum_{j \ne i} b_j}{N})(-\frac{N-1}{N})-\frac{N-1}{N}d_i && \nonumber\\
	+\frac{N-1}{N}(-\frac{b_i}{Na}-\frac{\sum_{j \ne i} b_j}{Na}) -\frac{1}{Na}(\frac{N-1}{N}b_i-\frac{\sum_{j \ne i}b_j}{N})& = & 0 \nonumber\\
	\eq
	
	With \eqref{eq:lambda} and \eqref{eq:p_i}, it is equivalent to
	\bsq
	\label{eq:Sharing-KKT2}
	\bq
	2c_ip_i+d_i-\lambda_c+\frac{1}{(N-1)a}(a\lambda_c+b_i) & = & 0 \label{eq:Sharing-KKT2.1}\\
	-\frac{b_i}{Na}-\frac{\sum \nolimits_{j \ne i} b_j}{Na} & = & \lambda_c \label{eq:Sharing-KKT2.2}\\
	D_i-\frac{N-1}{N}b_i+\frac{\sum \nolimits_{j \ne i} b_j}{N} & = & p_i \label{eq:Sharing-KKT2.3}
	\eq
	\esq
	
	Problem \eqref{eq:Central} is also a strictly convex optimization problem and the KKT condition is
	\bsq
	\label{KKT-central}
	\bq
	2\left[{c_i} - \frac{1}{{2(N-1)a}}\right]{p_i} + {d_i} + \frac{{{D_i}}}{{(N-1)a}} + \xi  = & 0 \label{KKT-central.1}\\
	\sum\limits_i {{p_i}} = & \sum\limits_i {{D_i}} \label{KKT-central.2}
	\eq
	\esq
	
	$\Rightarrow$ :  If $(p^*,b^*)$ is the NE of the sharing game \eqref{eq:SMK}, it satisfies the condition \eqref{eq:Sharing-KKT2}. 
	
	Sum up \eqref{eq:Sharing-KKT2.3} for all $i$, we have
	\bq
	\sum \limits_i p_i^* = \sum \limits_i D_i
	\eq
	which means \eqref{KKT-central.2} is satisfied.
	
	By $a \times$ \eqref{eq:Sharing-KKT2.2} + \eqref{eq:Sharing-KKT2.3}, we have
	\bq
	\label{KKT.1}
	D_i-p_i^* = a\lambda_c^*+b_i^*
	\eq
	Substitute \eqref{KKT.1} into \eqref{eq:Sharing-KKT2.1}, we have
	\bq
	\label{KKT.2}
	2[c_i-\frac{1}{2(N-1)a}]p_i^*+d_i+\frac{D_i}{(N-1)a}-\lambda_c^* =0
	\eq
	Let $\xi:=-\lambda_c^*$, constraint \eqref{KKT-central.1} is satisfied. As a result, $p^*$  is the optimal solution of \eqref{eq:Central}. 
	
	Obviously, we have $b_i^*=D_i-p_i^*-a\lambda_c^*$. $\lambda_c^*=\frac{1}{N}\sum \limits_i (2c_ip_i^*+d^i)$ will be proved latter in Appendix C.
	
	$\Leftarrow$: If $p^*$ is the optimal solution of problem \eqref{eq:Central}, then by letting
	\bsq
	\label{eq:construct}
	\bq
	p_i & = & p_i^* \\
	\lambda_c & = & -\xi^* = \frac{1}{N} \sum \limits_i (2c_ip_i^*+d_i) \label{eq:construct.1}\\
	b_i& = & D_i - p_i^* - a \lambda_c \label{eq:construct.2}
	\eq
	\esq 
	\eqref{eq:Sharing-KKT2.1} is obviously satisfied. Sum up \eqref{eq:construct.2} and together with \eqref{KKT-central.2} gives \eqref{eq:Sharing-KKT2.2} and \eqref{eq:Sharing-KKT2.3}.
	It satisfied the condition \eqref{eq:Sharing-KKT2}. This completes the proof.
\end{proof}

\renewcommand{\theequation}{B.\arabic{equation}}
\renewcommand{\thetheorem}{B.\arabic{theorem}}
\section{Proof of Proposition \ref{Thm:prop-3}}
\begin{proof}
	Under SDF-based sharing mechanism, given other prosumers’ strategies ${b_j},j \ne i$, by choosing  
	\bq
	{p_i} = {D_i}, \; {b_i} = \frac{{\sum\nolimits_{j \ne i} {b_j} }}{{(N - 1)}} \nonumber
	\eq
	with
	\bq
	{\lambda _c} =  - \frac{{\sum\nolimits_{j \ne i} {b_j} }}{(N - 1)a} \nonumber
	\eq
	We have $\Pi_i(p_i,b)=f_i(D_i)$, which means prosumer $i$ can acheive the same cost as under individual decision-making. Because each prosumer solves a minimization problem, so that we always have $\Pi(p_i^*,b^*) \le f_i(D_i)$. In consequence, a Pareto improvement is achieved for all prosumers.  
	
	If $p_i^*=D_i$ does not hold for all $i$, then as $p^*$ is the unique optimal solution of problem \eqref{eq:Central}, we always have
	\bq
	\sum \limits_i (c_i (p_i^*)^2+ d_i p_i^*) -\sum \limits_i \frac{(p_i^*-D_i)^2}{2(N-1)a} < \sum \limits_i (c_i D_i^2 +d_i D_i)
	\eq
	and becasue $\sum \limits_i \Pi(p_i^*, b^*) = \sum \limits_i f_i(p_i^*)$, so we have
	\bsq
	\bq
	&& \sum \limits_i \Pi(p_i^*,b^*) \nonumber\\
	<  && \sum \limits_i f_i(D_i) + \sum \limits_i \frac{(p_i^*-D_i)^2}{2(N-1)a} \nonumber \\
	< && \sum \limits_i f_i(D_i) 
	\eq
	\esq
	so that at least one strict inequality of \eqref{eq:pareto} holds. This completes the proof.
\end{proof}

\renewcommand{\theequation}{C.\arabic{equation}}
\renewcommand{\thetheorem}{C.\arabic{theorem}}

\section{Proof of Proposition \ref{Thm:prop-1}}
\begin{proof}
	If $(p^*,b^*)$ is the NE of the sharing game \eqref{eq:SMK}, then it satisfies the conditions \eqref{eq:Sharing-KKT2}.   
	
	1) Sum up \eqref{eq:Sharing-KKT2.3} for all $i$, and together with \eqref{KKT.1}, we have
	\bq
	\label{eq:sumq}
	\sum \limits_i (a\lambda_c^*+b_i^*) = 0
	\eq
	Sum up \eqref{eq:Sharing-KKT2.1} for all $i$ and with \eqref{eq:sumq}, we have
	\bq
	\lambda_c^*= \frac{1}{N} \sum\limits_i(2c_i p_i^* + d_i) = \frac{1}{N} \sum \limits_i md_i(p_i^*)
	\eq
	
	2) Given any $b_j, j \ne i$, let $b_i(p_i)=\frac{-Np_i+\sum \limits_{j \ne i}b_j+ND_i}{N-1}$ be the unique solution of constaints \eqref{SMK.2}, \eqref{MCC}. Then, for any fixed $b_j, j \ne i$, $\Pi_i(p_i,b_i,b_{-i})=\Pi_i(p_i,b_i(p_i),b_{-i})$.
	
	Derivatives of \eqref{SMK.2} and \eqref{MCC}  with respect to $p_i$ are
	\bq
	1+a \frac{\partial \lambda_c}{\partial p_i} + \frac{\partial b_i}{\partial p_i}=0
	\eq
	\bq
	Na \frac{\partial \lambda_c}{\partial p_i}+ \frac{\partial b_i}{\partial p_i}=0
	\eq
	Solving the equations we have
	\bq
	\frac{\partial \lambda_c}{\partial p_i} & = & \frac{1}{(N-1)a} \nonumber
	\eq
	\bq
	\frac{\partial b_i}{\partial p_i} &= & -\frac{N}{N-1} \nonumber
	\eq
	In consequence, the derivative of $\Pi_i(p_i,b_i(p_i),b_{-i})$ is 
	\bq
	\label{eq:derivative}
	\!\!&&\!\!\!\frac{\partial \Pi_i}{\partial p_i}(p_i,b_{-i}) \nonumber\\
	= \!\!&&\!\!\! (2c_i p_i+d_i)-\left[\lambda_c(b_i(p_i),b_{-i})-\frac{a \lambda_c(b_i(p_i),b_{-i})+b_i(p_i)}{(N-1)a}\right] \nonumber\\
	\eq
	The first term $2c_ip_i+d_i$ is the marginal disutility of prosumer $i$; the second term $\lambda_c$ is the marginal cost he needs to pay when buying from the market regardless of the mutual impact of different prosumers; the third term $(a \lambda_c+b_i)/{a(N-1)}$ reflects the marginal profit deviation due to the interest conflicts among different prosumers.
	
	If $md_i(p_i^*)>\lambda_c(b^*)$ and $q_i(b^*)<0$, then prosumer $i$ can always choose $\Delta p_i<0$ to reduce his cost, and thus, it is not a stable situation. Similarly, if $md_i(p_i^*)<\lambda_c(b^*)$ and $q_i(b^*)>0$, the market is also unstable as prosumer $i$ always has a better choice of $\Delta p_i>0$. As a result, $md_i(p_i^*)>\lambda_c(b^*)$ if and only if $q_i(b^*)>0$.  This completes the proof.
\end{proof}

\renewcommand{\theequation}{D.\arabic{equation}}
\renewcommand{\thetheorem}{D.\arabic{theorem}}
\section{Proof of Proposition \ref{Thm:prop-4}}
\begin{proof}
	1) According to Proposition \ref{Thm:prop-2}, $p_i^*$ is the unique solution of \eqref{eq:Central} and satisfies $\sum \limits_i p_i^* = \sum \limits_i D_i$. So $p_i^*$ must be a feasible solution to the social optimal problem \eqref{eq:AGG}. As $\bar p_i$ is the optimal solution of \eqref{eq:AGG}, we always have
	$$\sum \limits_{i \in \mathcal{I}} f_i(p_i^*) \ge \sum \limits_{i \in \mathcal{I}} f_i(\bar p_i)$$
	2) With the KKT conditions, we can obtain the optimal solutions to problem \eqref{eq:Central} and \eqref{eq:AGG} are
	\bqn
	p_i^*=\left(\sum \limits_j D_j +\sum \limits_j \frac{d_j+\frac{D_j}{a(N-1)}}{2c_j-\frac{1}{a(N-1)}}\right)\frac{\frac{1}{2c_i-\frac{1}{a(N-1)}}}{\sum \limits_j \frac{1}{2c_j-\frac{1}{a(N-1)}}}-\frac{d_i+\frac{D_i}{a(N-1)}}{2c_i-\frac{1}{a(N-1)}}
	\eqn
	and
	\bqn
	\bar p_i=\left(\sum \limits_j D_j+\sum \limits_j \frac{d_j}{2c_j}\right)\frac{\frac{1}{2c_i}}{\sum \limits_j \frac{1}{2c_j}}-\frac{d_i}{2c_i}
	\eqn
	respectively. Moreover, we have 
	\bsq
	\bq
	\bar p_i & \le & \left(N \bar D + N \frac{\bar d}{2 \underline c}\right) \frac{\bar c}{\underline c N} - \frac{\underline d}{2\bar c} \nonumber\\
	& = & (\bar D + \frac{\bar d}{2 \underline c})\frac{\bar c}{\underline c} - \frac{\underline d}{2 \bar c} \\
	\bar p_i & \ge & (N \underline D + N \frac{\underline d}{2\bar c})\frac{\underline c}{\bar c N}-\frac{\bar d}{2 \underline c} \nonumber\\
	& = & (\underline D+\frac{\underline d}{2 \bar c})\frac{\underline c}{\bar c}-\frac{\bar d}{2 \underline c}
	\eq
	\esq
	Letting $\bar A:= \sum \limits_j \frac{1}{2c_j}$ and $A^*:= \sum \limits_j \frac{1}{2c_j-\frac{1}{a(N-1)}}$ gives
	\bq
	0 \le \bar A- A^* & = & \sum \limits_j \frac{-\frac{1}{a(N-1)}}{2c_j(2c_j-\frac{1}{a(N-1)})} \nonumber\\
	& = & \sum \limits_j \frac{1}{2c_j-4ac_j^2(N-1)} \nonumber\\
	& \le & \frac{N}{2 \underline c -4 a\underline c^2(N-1)} \nonumber\\
	& \le & -\frac{1}{2 a\underline c^2}
	\eq
	Consequently, we have
	\bsq
	\bq
	\frac{\bar A- A^*}{\bar A} & \le & -\frac{\bar c}{a \underline c^2 N} \\
	\frac{\bar A - A^*}{\bar A} & \ge & 0
	\eq
	\esq
	Let $\bar B:=\sum \limits_j \frac{d_j}{2 c_j}$ and $B^*:=\sum \limits_j \frac{d_j + \frac{D_j}{a(N-1)}}{2c_j-\frac{1}{a(N-1)}}$. Then there are
	\bsq
	\bq
	\bar B- B^* & = & \sum \limits_j \frac{d_j+2c_jD_j}{-4ac_j^2(N-1)+2c_j} \nonumber\\
	& \le & \sum \limits_j \frac{|\bar d+2\bar {cD}|}{-4a\underline c^2(N-1)+2\underline c} \nonumber\\
	& \le & \frac{|\bar d+2\bar{cD}|}{-2 a\underline c^2}\\
	\bar B-B^* & = & \sum \limits_j \frac{d_j+2c_jD_j}{-4ac_j^2(N-1)+2c_j} \nonumber\\
	& \ge & \frac{-N|\underline d + 2\underline{cD}|}{-4a\underline c^2(N-1)+2\underline c} \nonumber\\
	& \ge & \frac{-|\underline d+ 2\underline{cD}|}{-2 a\underline c^2}
	\eq
	\esq
	Furthermore,  it is easy to see
	\bq
	|B^*| & \le & \sum \limits_j \left|\frac{d_j+\frac{D_j}{a(N-1)}}{2c_i-\frac{1}{a(N-1)}}\right| \nonumber\\
	& \le & \sum \limits_j \frac{|d_j|+|\frac{D_j}{a(N-1)}|}{2\underline c} \nonumber\\
	& \le & \sum \limits_j \frac{\bar d - max\{|\bar D|, |\underline D|\}/a}{2 \underline c} \nonumber\\
	& = & \frac{\bar d - max\{|\bar D|, |\underline D|\}/a}{2 \underline c}\cdot N
	\eq
	Let $\bar C:= \frac{d_i}{2c_i}$ and $C^* :=\frac{d_i+\frac{D_i}{a(N-1)}}{2c_i-\frac{1}{a(N-1)}}$. Direct calculation gives 
	\bsq
	\bq
	{\bar C- C^*} & =& \frac{d_i+2c_iD_i}{-4ac_i^2(N-1)+2c_i} \nonumber\\
	& \le & \frac{|\bar d+2\overline{cD}|}{-2 \underline c^2 aN} \\
	\bar C -C^* & = & \frac{d_i+2c_iD_i}{-4ac_i^2(N-1)+2c_i} \nonumber\\
	& \ge & \frac{-|\underline d+ 2 \underline{cD}|}{-4a\underline c^2 (N-1)+2\underline c} \nonumber\\
	& \ge & \frac{-|\underline d+ 2 \underline{cD}|}{-2 \underline c^2 aN} 
	\eq
	\esq
	Let $\bar E:=\frac{\frac{1}{2c_i}}{\bar A}$ and $E ^*: = \frac{\frac{1}{2c_i-\frac{1}{a(N-1)}}}{A^*}$, then we have
	\bsq
	\bq
	E^* & \le\frac{1}{2 \underline c A^*} & \le  \frac{2\bar c -1/a}{2\underline c} \cdot\frac{1}{N} \\
	E^* & \ge \frac{\frac{1}{2\bar c-1/a}}{A^*} & \ge \frac{2 \underline c}{2 \bar c-1/a} \cdot \frac{1}{N}
	\eq
	\esq 
	and
	\bsq
	\bq
	\frac{\bar E}{E^*} -1 & =&  \frac{A^*}{\bar A} \frac{2c_i}{2c_i -\frac{1}{a(N-1)}}-1 \nonumber\\
	& \le & (1- \frac{1}{(-2a \bar c+1)N})-1 \nonumber\\
	& \le &  - \frac{1}{\left(-2a\bar c+1\right) N}\\
	\frac{\bar E}{E^*} -1 & =&  \frac{A^*}{\bar A} \frac{2c_i}{2c_i -\frac{1}{a(N-1)}}-1 \nonumber\\
	& \ge & \left(1+\frac{\bar c}{a \underline c^2 N}\right) \left(1- \frac{1}{-a \underline cN}\right)-1 \nonumber\\
	& \ge &  \frac{\bar c}{a \underline c^2 N} + \frac{1}{a\underline c N}
	\eq
	\esq
	Therefore, there must be
	\bsq
	\bq
	\bar E - E^*  & \le &  \left( - \frac{1}{(-2a\bar c+1) N}\right) \frac{2\bar c-1/a}{2\underline c}\cdot\frac{1}{N} \\
	\bar E - E^*  & \ge & \left(\frac{\bar c}{a \underline c^2 N}+ \frac{1}{a\underline c N}\right)  \frac{2 \underline c}{2 \bar c-1/a} \cdot \frac{1}{N}
	\eq
	\esq
	The optimal solution of \eqref{eq:Central} and \eqref{eq:AGG} can be represented as
	$$\bar p_i =\left(\sum \limits_j D_j + \bar B\right) \bar E -\bar C$$
	$$p_i^*=\left(\sum \limits_j D_j + B^*\right) E^* -C^*$$
	First, we have
	\bq
	\left|\left(\sum \limits_j D_j +B^*\right)\left(\bar E- E^*\right)\right| & \le & \left|\sum \limits_j D_j +B^*\right| \left|\bar E- E^*\right| \nonumber\\
	& \le & \left|\sum \limits_j D_j + B^*\right| \frac{M_3}{N^2} \nonumber\\
	& \le & (\left|\sum \limits_j D_j\right|+\left|B^*\right|) \frac{M_3}{N^2} \nonumber\\
	& \le & (M_1 N +M_2 N) \frac{M_3}{N^2} \nonumber\\
	& = & \frac{M_3(M_1+M_2)}{N}
	\eq
	where 
	\bq
	M_1 & = & max\{|\underline D|,|\bar D|\} \nonumber\\
	M_2 & = & \frac{\bar d - max\{|\bar D|, |\underline D|\}/a}{2 \underline c}\nonumber\\
	M_3 & = & max\left\{\left|(\frac{1}{-2a\bar c+1})\frac{2\bar c-1/a)}{2\underline c}\right|,\right.\nonumber\\ && \left. \left|(\frac{\bar c}{ a \underline c^2}+ \frac{1}{a\underline c})\frac{2\underline c}{2\bar c-1/a}\right|\right \}
	\eq
	The difference between $\bar p_i$ and $p_i^*$ is
	\bq
	&& |\bar p_i - p_i^*| \nonumber\\
	& = & \left|\left(\sum \limits_j D_j +\bar B\right) \bar E -\bar C -\left(\sum \limits_j D_j + B^*\right) E^* +C^* \right|\nonumber\\
	& \le &   \left|\left(\sum \limits_j D_j +\bar B\right) \bar E - \left(\sum \limits_j D_j +  B^*\right) \bar E\right| \nonumber\\
	&& + \left|\left(\sum \limits_j D_j+  B^*\right) \bar E - \left(\sum \limits_j D_j + B^*\right) \bar E \right|+\left| \bar C - C^* \right|\nonumber\\
	& = & \left|\bar E (\bar B- B^*)\right| + \left|\left(\sum \limits_j D_j +B^*\right) \left(\bar E-E^*\right)\right| + \left|\bar C- C^*\right| \nonumber \\
	& \le & \frac{\bar c}{N \underline c} max\left\{\frac{|\bar d+2\bar{cD}|}{-2 a\underline c^2},\frac{-|\underline d+ 2\underline{cD}|}{-2 a\underline c^2}\right\} + M_3(M_1+M_2)\frac{1}{N} \nonumber\\
	&& + max\left\{\frac{|\bar d+2\overline{cD}|}{-2 \underline c^2 a}, \frac{|\underline d+ 2 \underline{cD}|}{-2 \underline c^2 a} \right\} \frac{1}{N}
	\eq 
	
	It is easy to see that  $|\bar p_i - p_i^*| \le \frac{\alpha}{N} $ holds for a large enough positive number $\alpha$.
	Because $\underline P \le \bar p_i \le \bar P$ and $-\frac{\alpha}{N} \le p_i^* -\bar p_i \le \frac{\alpha}{N}$, where
	$$\underline P=  \left(\underline D+\frac{\underline d}{2 \bar c}\right)\frac{\underline c}{\bar c}-\frac{\bar d}{2 \underline c}$$
	$$\bar P= \left(\bar D + \frac{\bar d}{2 \underline c}\right)\frac{\bar c}{\underline c} - \frac{\underline d}{2 \bar c}$$
	as a result, we have
	\bq
	2\underline P -\alpha \le 2\underline P- \frac{\alpha}{N} \le  p_i^*+\bar p_i \le 2 \bar P+\frac{\alpha}{N} \le 2\bar P+\alpha
	\eq
	
	For a given $\varepsilon>0$, we choose  a large enough number  $N_0:=\frac{1}{(\bar c \alpha max\{|2\underline P-\alpha|,|2\bar P+\alpha|\}+\alpha \bar d)\varepsilon}$. Then  for arbitrary number $N > N_0$, there is
	\bsq
	\bq
	&& \frac{1}{N}\left|\sum \limits_{i \in \mathcal{N}} f_i(p_i^*)-\sum \limits_{i \in \mathcal{N}} f_i(\bar p_i)\right| \nonumber\\
	& \le & \frac{1}{N} \sum \limits_i \left|c_i (p_i^*)^2 +d_i p_i -c_i(\bar p_i)^2 -d_i \bar p_i\right| \nonumber\\
	& \le & \frac{\bar c}{N} \sum \limits_i |p_i^*+\bar p_i||p_i^*-\bar p_i| + \frac{\bar d}{N} \sum \limits_i |p_i^*-\bar p_i| \nonumber\\
	& \le & (\bar c \alpha max\{|2\underline P-\alpha|,|2\bar P+\alpha|\}+\alpha \bar d)\frac{1}{N} \nonumber\\
	& < & \varepsilon
	\eq
	\esq
	This completes the proof.
	
\end{proof}

\renewcommand{\theequation}{E.\arabic{equation}}
\renewcommand{\thetheorem}{E.\arabic{theorem}}
\section{Proof of Proposition \ref{Thm:prop-4-2}}
\label{proof:prop-4-2}

\begin{proof}
	Suppose that $0<|a_1|<|a_2|$, and $p^{1*}$ corresponds to the NE of \eqref{eq:SMK} with $a=a_1$, and $p^{2*}$ corresponds to the NE with $a=a_2$. According to Proposition \ref{Thm:prop-2}, $p^{1*}$ and $p^{2*}$ are the unique optimal point of problem \eqref{eq:Central} under $a=a_1$ and $a=a_2$, respectively. Due to optimality,
	\bq
	&& \sum \nolimits_{i \in \mathcal{I}} f_i(p_i^{2*}) - \frac{\sum_{i \in \mathcal{I}}(p_i^{2*}-D_i)^2}{2(N-1)a_1} \nonumber\\
	& \ge & \sum \nolimits_{i \in \mathcal{I}} f_i(p_i^{1*}) - \frac{\sum_{i \in \mathcal{I}}(p_i^{1*}-D_i)^2}{2(N-1)a_1} 
	\eq
	which means
	\bq
	&& 2|a_1|(N-1)\left[\sum \nolimits_{i \in \mathcal{I}} f_i(p_i^{2*}) -\sum \nolimits_{i \in \mathcal{I}} f_i(p_i^{1*})\right]\nonumber\\
	& \ge & \left[\sum_{i \in \mathcal{I}}(p_i^{1*}-D_i)^2-\sum_{i \in \mathcal{I}}(p_i^{2*}-D_i)^2\right]
	\eq
	If we have
	$$
	\sum \nolimits_{i \in \mathcal{I}} f_i(p_i^{1*}) < \sum \nolimits_{i \in \mathcal{I}} f_i(p_i^{2*}) 
	$$
	then
	\bq
	&& 2|a_2|(N-1)\left[\sum \nolimits_{i \in \mathcal{I}} f_i(p_i^{2*}) -\sum \nolimits_{i \in \mathcal{I}} f_i(p_i^{1*})\right] \nonumber\\
	& \ge & 2|a_1|(N-1)\left[\sum \nolimits_{i \in \mathcal{I}} f_i(p_i^{2*}) -\sum \nolimits_{i \in \mathcal{I}} f_i(p_i^{1*})\right]\nonumber\\
	& \ge & \left[\sum_{i \in \mathcal{I}}(p_i^{1*}-D_i)^2-\sum_{i \in \mathcal{I}}(p_i^{2*}-D_i)^2\right]
	\eq
	which means
	\bq
	&& \sum \nolimits_{i \in \mathcal{I}} f_i(p_i^{2*}) - \frac{\sum_{i \in \mathcal{I}}(p_i^{2*}-D_i)^2}{2(N-1)a_2} \nonumber\\
	& \ge & \sum \nolimits_{i \in \mathcal{I}} f_i(p_i^{1*}) - \frac{\sum_{i \in \mathcal{I}}(p_i^{1*}-D_i)^2}{2(N-1)a_2} 
	\eq
	and is contradict to the assumption that $p^{2*}$ corresponds to the NE under $a=a_2$ (also the optimal solution of problem \eqref{eq:Central} under $a=a_2$), which completes the proof.
\end{proof}

\renewcommand{\theequation}{F.\arabic{equation}}
\renewcommand{\thetheorem}{F.\arabic{theorem}}
\section{Proof of Proposition \ref{Thm:prop-6}}
\label{proof:prop-6}
\begin{proof}
Given other prosumers' strategies, prosumer $i$ solves problem \eqref{eq:SMK-Ex} to obtain his optimal strategy. Use $b_i$ to represent $\lambda_c$, problem \eqref{eq:SMK-Ex} can be reduced to the following optimization problem.
\bsq
\label{eq:SMK-Ex1}
\bq
\mathop {\min }\limits_{{p_i^k,\forall k \in \mathcal{K}_i},{b_i}} && \!\!\!\!\! \sum\nolimits_k \left[{c_i^k}\left(p_i^k\right)^2 + {d_i^k}{p_i^k}\right] \nonumber\\
&& \!\!\!\!\!+ \left(\frac{I-1}{I}b_i-\frac{\sum \limits_{j \ne i}b_j}{I}\right)\left({\frac{-b_i-\sum \limits_{j \ne i} b_j}{Ia}}\right)\\
{\rm{s}}.{\rm{t}}.&&  \!\!\!\!\!\sum\nolimits_k {p_i^k} + \frac{I-1}{I}b_i - \frac{1}{I} \sum \limits_{j \ne i} b_j= {D_i} :\mu_i 
\eq
\esq
It is easy to check that problem \eqref{eq:SMK-Ex1} is a strictly convex optimization problem. Thus, the KKT condition below is the necessary and sufficient condition for the optimal point.
	\bsq
	\label{eq:coalition-KKT}
	\bq
	2c_i^kp_i^k + d_i^k + {\mu _i} & = & 0,\forall k  \label{eq:coalition-KKT.1}\\
	\frac{I-1}{I}\left(\frac{-b_i-\sum \limits_{j \ne i}b_j}{Ia}\right)
	+\frac{I-1}{I}\mu_i & & \nonumber\\
	-\frac{1}{Ia}\left(\frac{I-1}{I}b_i-\frac{\sum \limits_{j \ne i}b_j}{I}\right) & = & 0 \label{eq:coalition-KKT.2}\\
	\sum\nolimits_k {p_i^k} + \frac{I-1}{I}b_i - \frac{1}{I} \sum \limits_{j \ne i} b_j &= & {D_i} \label{eq:coalition-KKT.3}
	\eq
	\esq
	
	The KKT condition of problem \eqref{eq:Central-Ex} is
	\bsq
	\label{eq:Central-Ex-KKT}
	\bq
	\left[2c_i^k - \frac{1}{{(I - 1)a}}\right]p_i^k + d_i^k & &\nonumber\\
	+ \frac{{{D_i}}}{{(I - 1)a}} 
	- \frac{{\sum\limits_{j \in {K_i},j \ne k} {p_i^j} }}{{(I - 1)a}} & = & -\xi^{'} \label{eq:Central-Ex-KKT.1}\\
	\sum\limits_i \sum\limits_k {{p_i^k}} & = & \sum\limits_i {{D_i}} \label{eq:Central-Ex-KKT.2}
	\eq
	\esq
	$\Rightarrow$: If $(p^*,b^*)$ is the NE of the sharing problem for MRP \eqref{eq:SMK-Ex}, then it satisfies the KKT conditions \eqref{eq:coalition-KKT}. Sum up the \eqref{eq:coalition-KKT.3} for all $i$, \eqref{eq:Central-Ex-KKT.2} is met. Substitute \eqref{eq:coalition-KKT.2}-\eqref{eq:coalition-KKT.3} into \eqref{eq:coalition-KKT.1}, we have
	\bq
	\label{eq:kkt2}
	 && \left[2c_i^k - \frac{1}{{(I - 1)a}}\right]p_i^{k*} + d_i^k + \frac{{{D_i}}}{{(I - 1)a}} - \frac{{\sum\limits_{j \in {K_i},j \ne k} {p_i^{j*}} }}{{(I - 1)a}} \nonumber\\
	 & = & -\frac{1}{Ia}(b_i^*+\sum \limits_{j \ne i} b_j^*)
	\eq
	which is a constant for $\forall i, \forall k \in \mathcal{K}_i$. Let $\xi^{'}:= \frac{1}{Ia}(b_i^*+\sum \limits_{j \ne i} b_j^*)$, and \eqref{eq:Central-Ex-KKT.1} is met. In consequence, $p^*$ is also the optimal solution of problem \eqref{eq:Central-Ex}.
	
	Let $\lambda_c:= -\frac{1}{Ia}(b_i^*+\sum \limits_{j \ne i} b_j)$. Sum up \eqref{eq:kkt2} for $i$, we have
		$$\lambda_c^*(p^*)=\frac{1}{I}\sum \limits_i md_i (p_i^*)$$
	With \eqref{eq:coalition-KKT.2}, we have 
		$$b_i^*(p^*)=D_i-\sum \limits_k p_i^{k*}-a\lambda_c^*(p^*) $$
	$\Leftarrow$: if $p^*$ is the optimal solution of problem \eqref{eq:Central-Ex}, then it satisfies the KKT conditions \eqref{eq:Central-Ex-KKT}. If we let
	\bsq
	\label{eq:coalition-3}
	\bq
	\lambda_c & = & -\xi^{'*} = \frac{1}{I} \sum \limits_i {md}_i(p_i^*) \label{eq:coalition-3.1}\\
	\mu_i & = & \frac{1}{(I-1)a}(D_i-\sum_k p_i^{k*}) - \lambda_c \label{eq:coalition-3.2}\\
	b_i & = & D_i-(a\lambda_c+\sum\limits_k (p_i^k)^*) \label{eq:coalition-3.3} \\
	p_i & = & p_i^*
	\eq
	\esq
	
	Obviously, \eqref{eq:coalition-KKT.1} is satisfied. Sum up \eqref{eq:coalition-3.3} for all $i$ together with \eqref{eq:Central-Ex-KKT.2} gives that
	$$\lambda_c = -\frac{1}{Ia} \sum_i b_i$$
	It is easy to check \eqref{eq:coalition-KKT.2} and \eqref{eq:coalition-KKT.3} are satisfied. 
	Thus, $(p^*,b^*)$ is the NE of the sharing problem \eqref{eq:SMK-Ex}, which completes the proof.
\end{proof}

\renewcommand{\theequation}{G.\arabic{equation}}
\renewcommand{\thetheorem}{G.\arabic{theorem}}
\section{Proof of Proposition \ref{Thm:prop-8}}
\begin{proof}
	
	We consider a special situation, where all the $c_i,\forall i$ equal to a same value $c$, and each prosumer possess the same number of resource, which is $K_I$. $IK_I$ is a fixed value equals to $N$.
	
	With the KKT condition \eqref{eq:Central-Ex-KKT}, we have
	\bsq
	\label{eq:simplified-1}
	\bq
	2cp_i^k+d_i^k+\frac{D_i-\sum \limits_{k \in \mathcal{K}_i}p_i^k}{(I-1)a} & = & -\xi^{'} \label{eq:simplified-1.1} \\
	\sum \limits_i \sum \limits_k p_i^k & = & \sum \limits_i D_i \label{eq:simplified-1.2}
	\eq
	\esq
	Denote $2cp_i^k+d_i^k$ as $md_i^k$. It is obvious that for all $k \in \mathcal{K}_i$, $md_i^k$ are equal, and so we use $md_i$ to represent $md_i^k$ for all $k$.
	Sum up \eqref{eq:simplified-1.1} for all $i$ and $k$, we have
	\bq
	2c\sum\limits_i {\sum\limits_k {p_i^k} }  + \sum\limits_i {\sum\limits_k {d_i^k} }  + \frac{{K_I(\sum\limits_i {{D_i}}  - \sum\limits_i {\sum\limits_{k \in {\mathcal{K}_i}} {p_i^k} } )}}{{(I - 1)a}} + N\xi^{'}  = 0
	\eq
	Together with \eqref{eq:simplified-1.2}, it is easy to find that $\xi^{'}$ is independent of $I$ and $\sum \limits_i md_i + I \xi^{'} =0$. Assume that the optimal marginal disutility is $md_i^*$, then according to \eqref{eq:simplified-1.1}, we have
	\bq
	\frac{{{D_l} - \sum\limits_{k \in {\mathcal{K}_i}} {p_i^{k*}} }}{{(I - 1)a}} =  - md_i^* - \xi ^{'*}
	\eq
	The objective function \eqref{Central-Ex.1} can be rewritten as
	\bq
	\pi:= \pi_1+\pi_2  \nonumber
	\eq
	where
	\bq
	\pi_1  =  \sum\limits_i {\sum\limits_k {[c{{(p_i^{k*})}^2} + d_i^kp_i^{k*}]} } =  \sum\limits_i {\sum\limits_k {\frac{{(md_i^{*})^2 - {{(d_i^k)}^2}}}{{4c}}} } \nonumber
	\eq
	\bq
	{\pi _2} = -\frac{{\sum\limits_i {{{({D_i} - \sum\limits_{k \in {\mathcal{K}_i}} {p_i^{k*}} )}^2}} }}{{2(I - 1)a}} = -\frac{{(I - 1)a}}{2}\sum \limits_i {({md_i^*} + \xi^{'*} )^2} \nonumber
	\eq
	Then we have
	\bq
	\label{eq:simplified-2}
	\pi_1 & = & \sum\limits_i {\sum\limits_k {\frac{{(md_i^*)^2 - {{(d_i^k)}^2}}}{{4c}}} } \nonumber\\
	& = & \frac{K_I}{4c} \sum \limits_i (md_i^*)^2 - \sum \limits_i \sum \limits_k \frac{(d_i^k)^2}{4c} \nonumber \\
	& = & \frac{K_I}{4c} \sum \limits_i (md_i^*+\xi^{'*})^2 + \frac{N(\xi^{'*})^2-\sum \limits_i \sum \limits_k (d_i^k)^2}{4c}
	\eq
	The first term of \eqref{eq:simplified-2} is variational and the second term is a constant. So the change of $\pi_1$ is related with $\frac{K_I}{4c} \sum \limits_i (md_i^*+\xi^{'*})^2 $, the change of $\pi_2$ is related with $ -\frac{{(I - 1)a}}{2} \sum \limits_i {({md_i^*} + \xi^{'*} )^2}$ and the change of $\pi$ is related with $\frac{K_I-2ac(I-1)}{4c} \sum \limits_i (md_i^*+\xi^{'*})^2$. Obviously, $\frac{K_I-2ac(I-1)}{4c}$ is always positive for all $I>1$.
	
	Next, we define an equal partition such that $\pi$ is decreasing. Assume that there is $I$ multi-resource prosumers, and the optimal output is $p_i^{k*}$. Then we introduce competition and allocate the resources owned by one prosumer and its demand to $Z \in \mathbb{Z^{+}}$  prosumers such that $I^{'}=ZI$ and $K_I=ZK_{I^{'}}$, and satisfies $\forall i, \forall z_1,z_2 \in \{1,...,Z\}$
	\bq
	({D_{Z(i-1)+z_1}^{'}} - \sum\limits_{k=1}^{K_{I^{'}}} {p_{_{Z(i-1)+z_1}}^{k*}} )({D_{Z(i-1)+z_2}^{'}} - \sum\limits_{k=1}^{K_{I^{'}}} {p_{_{Z(i-1)+z_2}}^{k*}} ) & \ge & 0 \nonumber\\
	\sum \limits_{z=1}^Z D_{{Z(i-1)+z}}^{'} &= & D_i \nonumber
	\eq
	Then
	\bsq
	\label{eq:simplified-3}
	\bq
	\pi^{I*} & = & \pi_1^{I*} - \sum\limits_i {\frac{{{{({D_i} - \sum\limits_{k=1}^{K_I} {p_i^{k*}} )}^2}}}{{(I - 1)a}}} \nonumber\\
	& = &  \pi_1^{I*}  - \sum \limits_i \frac{(\sum \limits_{z=1}^Z D_{Z(i-1)+z}^{'}- \sum \limits_{z=1}^Z \sum \limits_{k=1}^{K_{I^{'}}} p_{Z(i-1)+z}^{k*})^2}{(I-1)a} \nonumber \\
	& \ge & \pi_1^{I*} - \sum \limits_i \frac{\sum \limits_{z=1}^Z (D_{Z(i-1)+z}^{'}-\sum \limits_{k=1}^{K_{I^{'}}}p_{Z(i-1)+z}^{k*})^2}{(I-1)a} \nonumber\\
	& \ge & \pi_1^{I*} - \sum \limits_i \frac{\sum \limits_{z=1}^Z (D_{Z(i-1)+z}^{'}-\sum \limits_{k=1}^{K_{I^{'}}} p_{Z(i-1)+z}^{k*})^2}{(I^{'}-1)a} \nonumber\\
	& \ge & \pi^{I^{'}*} 
	\eq
	\esq
	
	\eqref{eq:simplified-3} tells us that the objective function $\pi$ decreases with $I$. When $I$ grows, the coefficient of $\pi_1$ (which is $\frac{K_I}{4c}$) decreases while the coefficient of $\pi_2$ (which is $\frac{-(I-1)a}{2}$) increases. If $\sum \limits_i (md_i+\xi^{'})^2$ decreases with $I$, then obviously $\pi_1$ decreases with $I$. Otherwise, if $\sum \limits_i (md_i+\xi^{'})^2$ increases with $I$, $\pi_2$ increases with $I$. If $\pi_1$ also increases, then $\pi^{I*}<\pi^{I^{'}*}$, which is contradict to \eqref{eq:simplified-3}. In conclusion, we always have the social total cost $\pi_1$ reduces with $I$, which means introducing effective competition can improve social welfare. 
	
	Moreover, as $\frac{K-2ac(I-1)}{4c}$ is always positive and reaches the minimum when $I=\sqrt{{-2ac}/{N}}$. If $\sqrt{{-2ac}/{N}} \le 1$, then $\frac{K-2ac(I-1)}{4c}$ increases with $I>1$, which tells that $\sum \limits_i (md_i+\xi^{'})^2$ is decreasing in $I$, so that
	\bq
	Var(md_i^{k*},I) & = & \frac{1}{N}\sum \limits_{i=1}^I \sum \limits_{k=1}^{K_I} (md_i^{k*}+\xi^{'})^2 \nonumber\\
	& = & \frac{K_I}{N} \sum \limits_{i=1}^{I} (md_i+\xi^{'})^2 \nonumber\\
	& \ge &  \frac{K_{I^{'}}}{N} \sum \limits_{i=1}^{I^{'}} (md_i+\xi^{'})^2= Var(md_i^{k*}, I^{'}) \nonumber
	\eq
	This completes the proof.
\end{proof}

\end{document}